\documentclass[11pt]{article}
\usepackage{amsmath,amsfonts,amssymb}
\usepackage{bbm}
\usepackage{cases}
\usepackage[hidelinks]{hyperref}

\def\<{{\langle}}
\def\>{{\rangle}}
\allowdisplaybreaks[4]

\newtheorem{theorem}{Theorem}[section]
\newtheorem{lemma}[theorem]{Lemma}
\newtheorem{corollary}[theorem]{Corollary}

\newtheorem{remark}[theorem]{Remark}

\newtheorem{definition}[theorem]{Definition}

\begin{document}

\title{\vspace*{-1.5cm}
Large deviations of stochastic heat equations with logarithmic nonlinearity}

\author{Tianyi Pan$^{1}$, Shijie Shang$^{1}$,
Tusheng Zhang$^{2}$}
\footnotetext[1]{\, School of Mathematics, University of Science and Technology of China, Hefei, China. Email: pty0512@mail.ustc.edu.cn (Tianyi Pan), sjshang@ustc.edu.cn (Shijie Shang).}
\footnotetext[2]{\, School of Mathematics, University of Manchester, Oxford Road, Manchester M13 9PL, England, U.K. Email: tusheng.zhang@manchester.ac.uk}
\maketitle

\begin{abstract}
In this paper, we establish a large deviation principle for the solutions to the stochastic heat equations with logarithmic nonlinearity driven by Brownian motion, which is neither locally Lipschitz nor locally monotone.
Nonlinear versions of  Gronwall's inequalities and Log-Sobolev inequalities play an important role.
\end{abstract}


\medskip


\section{Introduction}

\quad In this paper, we study the small noise large deviation principle (LDP) of stochastic heat equations with logarithmic nonlinear drift term driven by Brownian motion, which is written as follows,
\begin{numcases}{}
 du(t,x) = \Delta u(t,x) dt+ u(t,x)\log|u(t,x)| dt+ \sigma(u(t,x))dW_t, \   t>0, x\in D , \nonumber\\
 u(t,x)=0, \quad t>0, x\in\partial D , \nonumber\\
\label{1.a} u(0,x)=u_0(x), \quad  x\in D ,
\end{numcases}
where $D$ is a bounded domain in $\mathbb{R}^d$ with smooth boundary. The coefficient $\sigma(\cdot): \mathbb{R}\rightarrow \mathbb{R}$ is a deterministic continuous function. $W$ is an 1-dimensional standard Brownian motion defined on a complete  filtrated probability space $(\Omega, {\cal F}, {\cal F}_t, P)$.  In this paper, we assume that the initial value $u_0$ is a deterministic element in $L^2(D)$. 

\vskip 0.3cm

We like to point out  that the drift coefficient $b(z)=z\log|z|$ neither has linear growth nor be locally Lipschitz. In fact, this function $b$ is locally $\log$-Lipschitz and  of superlinear growth. And equation (\ref{1.a}) does not fall into the category of stochastic partial differential equations with locally monotone coefficients known in the literature. We refer the readers to \cite{SZ} and references therein.

\vskip 0.3cm

PDEs with a logarithmic nonlinearity have been  introduced in the study of nonlinear wave mechanics. The logarithmic wave mechanics and logarithmic Schr\"{o}dinger equations have been studied by many authors, see \cite{R,BM}. The logarithmic deterministic parabolic equations have also been widely studied, we refer the readers to \cite{CLL, CT, JYC, DZ}  and references therein  for details.

\vskip 0.3cm

The stochastic heat equations with logarithmic nonlinearity driven by Brownian motion was studied in the paper \cite{SZ} by Shang and Zhang.
They proved that equation (\ref{1.a}) has a unique global strong probabilistic solution when the diffusion coefficient $\sigma$ satisfies a locally Lipschitz condition and a certain superlinear growth.
We would like also to mention the paper \cite{DKZ} where the authors studied the stochastic reaction diffusion equations on the interval $[0, 1]$ driven by space-time noise with coefficients which are locally Lipschitz and satisfy some superlinear logarithmic growth.
Stochastic reaction diffusion equations on the whole line $\mathbb{R}$ with logarithmic nonlinearity driven by space-time white noise were considered in \cite{SZ2}

\vskip 0.3cm

The large deviation theory has wide applications in many areas, e.g. statistical mechanics, risk management and hydrodynamics, see \cite{AO} for more applications. For the LDP of stochastic differential equations (SDEs) and stochastic partial differential equations (SPDEs), there exists a large amount of literature, we refer the readers to \cite{CR,FZ,L,SS,BDM,P,DXZZ,DWZZ,MSZ} and references therein for details.

\vskip 0.3cm
The purpose of this paper is to establish a Freidlin--Wentzell type LDP for the stochastic heat equation (\ref{1.a}). To obtain the LDP of the solutions, we will adopt the weak convergence method introduced in \cite{BD}. Especially, we will use the more convenient sufficient conditions given in the paper \cite{MSZ}, see Theorem \ref{criterion} below.  To this end, we first establish the well-posedness of the corresponding skeleton equation using the Galerkin approximations, where two versions of nonlinear Gronwall's inequality played an important role.
To verify the  conditions given in \cite{MSZ}, we first show they hold on a small time interval again using a version of nonlinear Gronwall's inequality and then extend it to the whole interval $[0,T]$  by induction. For the weak convergence,  our first step is to prove it under a weaker metric and then we improve the convergence  under a stronger metric.

\vskip 0.3cm
The rest of the paper is organized as follows. In Section 2, we present the framework for (\ref{1.a}), give our hypotheses and introduce the weak convergence method of LDP. Section 3 is  to introduce the main result. In Section 4, we establish the well-posedness of the skeleton equation. In Section 5, we prove the large deviation result. Section 6 is the Appendix containing two versions of nonlinear Gronwall's inequalities used in this paper.

\section{Framework}\label{S:2}
\setcounter{equation}{0}

\quad In this section, we will set up the framework and briefly recall the weak convergence method  in the large deviation principle theory.
 Let $H:=L^2(D)$ be the $L^2$ space with the norm and the inner product respectively denoted by $\|\cdot\Vert$ and $(\cdot\ ,\ \cdot)$. Denote the Sobolev space $H_0^1(D)$ by $V$, which is the completion of the space of compactly supported smooth functions $C^{\infty}_c(D)$ under the norm
\begin{align}
  \big\|u\big\|_{V}^2=\int_D|\nabla u(x)|^2dx.
\end{align}
There exists an orthonormal basis $\{e_i\}_{i=1}^{\infty}$ of $H$ consisting of the eigenvectors of the negative Laplace operator under zero boundary conditions with the corresponding eigenvalues $0<\lambda_i \uparrow\infty$, satisfying:
\begin{align}\label{A eigenvalue}
\Delta e_i=-\lambda_i e_i, \quad e_i|_{\partial D}=0,\quad i\in\mathbb{N}.
\end{align}
Moreover, $\{e_i\}_{i=1}^{\infty}$ is an orthogonal basis of $V$.
Recall the $\rm Poincar\acute{e}$ inequality, i.e.
\begin{align}\label{Poincare}
\big\|u\big\|^2 \leq \frac{1}{\lambda_1}\big\|u\big\|_{V}^2, \quad \forall\, u\in V .
\end{align}
In order to deal with the logarithmic term, we need the logarithmic Sobolev inequality (see \cite{G})  of the following form.
For any $\varepsilon>0$ and $u\in V$, we have
\begin{align}\label{log sobolev}
\int_D |u(x)|^2\log|u(x)| dx \leq \varepsilon \big\|u \big\|_{V}^2 + \left(\frac{d}{4}\log\frac{1}{\varepsilon}\right) \big\|u\big\|^2 + \big\|u\big\|^2\log\big\|u\big\|.
\end{align}
Set
\[
\log_{+}z:=\log (1\vee z) .
\]
From the above logarithmic Sobolev inequality, it follows that for any $\varepsilon>0$ and $u\in V$,
\begin{align}\label{log sobolev modi}
& \int_D |u(x)|^2\log_{+}|u(x)| dx \nonumber\\
\leq & \varepsilon \big\|u \big\|_{V}^2 + \left(\frac{d}{4}\log\frac{1}{\varepsilon}\right) \big\|u\big\|^2 + \big\|u\big\|^2\log\big\|u\big\|+ \frac{1}{2\mathrm{e}}m(D) ,
\end{align}
where $m(D)$ denotes the Lebesgue measure of the domain $D$.

Identifying the Hilbert space $H$ with its dual $H^*$ via the Riesz representation, we obtain a Gelfand triple
\begin{align*}
V\subset H\subset V^*.
\end{align*}

\noindent Denote by $\langle f,v\rangle$ the canonical dual pairing between $f\in V^*$ and $v\in V$. We have
\begin{align}\label{eq P4 star}
(u,v)=\langle u,v\rangle,\ \ \ \forall\,u\in H,\ \ \forall\,v\in V.
\end{align}

\noindent Set
\begin{align*}
	u(t)(x):= & u(t,x), \quad\quad \left(u(t)\log|u(t)|\right)(x):= u(t,x)\log|u(t,x)|, \\
 \sigma(u(t))(x):= & \sigma(u(t,x)) .
\end{align*}
Then (\ref{1.a}) can be reformulated as the following stochastic evolution  equation
\begin{numcases}{}
  u(t) = u_0 +\int_0^t \Delta u(s)ds + \int_0^t u(s)\log|u(s)| ds + \int_0^t \sigma(u(s))dW_s , \nonumber\\
\label{Abstract}  u(0)=u_0\in H.
\end{numcases}

\begin{definition}
An $H$-valued $\{\mathcal{F}_t\}$-adapted continuous  stochastic process $(u_t)_{t\geq 0}$ is called a  solution of (\ref{Abstract}), if the following two conditions hold:
\begin{itemize}
  \item [(i)] $u\in L^2([0,T];V)$ for any $T>0$, $P$-a.s.
  \item [(ii)] $u$ satisfies the equation (\ref{Abstract}) in $V^*$, $P$-a.s for any $t\geq 0$.
\end{itemize}
\end{definition}

\vskip 0.3cm

Now we introduce the hypotheses on the diffusion coefficient $\sigma$.
\begin{itemize}
  \item [\hypertarget{H}{{\bf (H)}}] There exist positive constants $L_1$ and $L_2$  such that for all $x,y\in \mathbb{R}$,
\begin{gather}
\label{Local-Lipschitz} |\sigma(x)-\sigma(y)| \leq  L_1|x-y|+ L_2 |x-y| \left(\log_{+}(|x|\vee |y|)\right)^{\frac{1}{2}} .
\end{gather}
\end{itemize}
\begin{remark}
  The assumption (\hyperlink{H}{H}) implies that  $\sigma$ is continuous and there exist positive constants $L_3$, $L_4$  such that for all $x\in \mathbb{R}$,
\begin{equation}\label{H2}
|\sigma(x)| \leq  L_3 + L_4 |x|\left(\log_{+} |x|\right)^{\frac{1}{2}}.
\end{equation}
\end{remark}

The following result is the Theorem 6.6 in \cite{SZ} giving the well-posedness  of equation (\ref{Abstract}).
\begin{theorem}\label{2.2}
	Supposing hypothesis (\hyperlink{H}{H}) hold. Then there exists a unique global solution u to (\ref{Abstract}) for every initial value $u_0\in H$.
\end{theorem}

Next, we turn to the definition of LDP.

\begin{definition}
	Let $\mathcal{E}$ be a Polish space with the Borel $\sigma$-field $\mathcal{B}(\mathcal{E})$. A function I: $\mathcal{E}\rightarrow[0,+\infty]$ is called a \textbf{rate function} if I is lower semicontinuous and the level set $\{e\in \mathcal{E}:I(e)\leq M\}$ is a compact subset of $\mathcal{E}$ for each $M<\infty$.
\end{definition}

\begin{definition}
	A family of $\mathcal{E}$-valued random variables $\{X^\varepsilon\}_{\varepsilon>0}$ is said to satisfy the LDP on $\mathcal{E}$ with rate function I if for each Borel subset B of $\mathcal{E}$,
	$$-\inf\limits_{e\in\mathring{B}}I(e)\leq \liminf\limits_{\varepsilon\rightarrow 0}\varepsilon^2 \log P(X^\varepsilon\in B)\leq\limsup\limits_{\varepsilon\rightarrow 0}\varepsilon^2 \log P(X^\varepsilon\in B)\leq -\inf\limits_{e\in\bar{B}}I(e).$$
\end{definition}

 Next we will introduce a sufficient condition for the LDP for a sequence of Wiener functionals.

 Let $\{W_t\}_{t\geq 0}$ be a real-valued Wiener process on a complete filtrated probability space $(\Omega,\mathcal{F},\mathcal{F}_t,P)$. Suppose for each $\varepsilon>0$, $\mathcal{G}^\varepsilon:C([0,T];\mathbb{R})\rightarrow \mathcal{E}$ is a measurable map and $X^\varepsilon=\mathcal{G}^\varepsilon(W.)$. For $N>0$, we set
\begin{align*}
\mathcal{A}&=\Big\{v:v \text{ is a real-valued  $\{\mathcal{F}_t\}$-predictable process such that}\\& \quad\quad\quad
\int_{0}^{T}|v(s,\omega)|^2ds<\infty,\  P\text{-a.s}\Big\},\\
S_N&=\Big\{\phi\in L^2([0,T],\mathbb{R}):\int_{0}^{T}|\phi(s)|^2ds \leq N\Big\},\\
\mathcal{A}_N&=\Big\{v\in \mathcal{A}:v(\cdot,\omega)\in S_N,\  P\text{-a.s}\Big\}.
	\end{align*}
To establish the LDP of the Wiener functionals $\{X^\varepsilon\}_{\varepsilon>0}$, we will use the following sufficient conditions established in \cite{MSZ}, which are based on a criteria of Budhiraja-Dupuis in \cite{BD}.

\begin{theorem}{}\label{criterion}
	If there exists a measurable map $\mathcal{G}^0:C([0,T];\mathbb{R})\rightarrow \mathcal{E}$ such that the following two conditions hold\\
	(a) For every $N<\infty$, for any family $\{h_\varepsilon\}_{\varepsilon>0}\subseteq \mathcal{A}_N$ and for any $\delta>0$,$$\lim\limits_{\varepsilon\rightarrow 0}P(\rho(Y^{h_\varepsilon},X^{h_\varepsilon})>\delta)=0$$
	where $X^{h_\varepsilon}:=\mathcal{G}^\varepsilon(W.+\frac{1}{\varepsilon}\int_{0}^{.}h_\varepsilon(s)ds)$, $Y^{h_\varepsilon}:=\mathcal{G}^0(\int_{0}^{.}h_\varepsilon(s)ds)$ and $\rho(\cdot,\cdot)$ stands for the metric of the space $\mathcal{E}$.\\
	(b) For every $N<\infty$ and any family $\{h_\varepsilon\}_{\varepsilon > 0}\subseteq S_N$ that converges weakly to some element h in $L^2([0,T];\mathbb{R})$ as $\varepsilon \rightarrow 0$ ,  we have $Y^{h_\varepsilon}\rightarrow Y^{h}\ in\ \mathcal{E}$.\\
Then the family $X^\varepsilon=\mathcal{G}^\varepsilon(W_\cdot)$ satisfies a large deviation principle with the rate function
	\begin{align*}
	I(f)=\inf\limits_{\{v\in L^2([0,T];\mathbb{R}):f=\mathcal{G}^0(\int_{0}^{\cdot}v(s)ds)\}}\Big\{\frac{1}{2}\int_{0}^{T}|v(s)|^2ds\Big\},
	\end{align*}
 with the convention $\inf\{\emptyset\}=\infty$.
\end{theorem}

\section{Statement of the main result}
     \quad In the remainder of this paper, we will take the Polish space $$\mathcal{E}=C([0,T];H)\cap L^2([0,T];V)$$ equipped with the metric,
     $$\rho(u,v)^2=\sup\limits_{s\in[0,T]}\big\|u(s)-v(s)\big\|^2+\int_{0}^{T}\big\|u(s)-v(s)\big\|_V^2ds,\ \forall\ u,v\in \mathcal{E}.$$

     Assume that (\hyperlink{H}{H}) is satisfied, then from Theorem \ref{2.2}, there exists a unique solution to the following equation:
    \begin{numcases}{}\label{eps}
    	u^\varepsilon(t) = u_0 +\int_0^t \Delta u^\varepsilon(s)ds + \int_0^t u^\varepsilon(s)\log|u^\varepsilon(s)| ds + \varepsilon\int_0^t \sigma(u^\varepsilon(s))dW_s , \nonumber\\
    	\label{Abstracte}  u^\varepsilon(0)=u_0\in H.
    \end{numcases}
\indent
    By the Yamada-Watanabe theorem in $\cite{RSZ}$, the solution of (\ref{eps}) determines a measurable map $\mathcal{G}^\varepsilon:C([0,T];\mathbb{R})\rightarrow \mathcal{E}$ such that for all standard Brownian motion $W$, $\mathcal{G}^\varepsilon(W)$ is the unique solution of (\ref{eps}).\\
\indent
    We also need to consider the so-called skeleton equation
    \begin{numcases}{}\label{Skeleton}
    	u^h(t) = u_0 +\int_0^t \Delta u^h(s)ds + \int_0^t u^h(s)\log|u^h(s)| ds + \int_0^t h(s)\sigma(u^h(s))ds , \nonumber\\
    	 u^h(0)=u_0\in H,
    \end{numcases}
where $h\in L^2([0,T];\mathbb{R})$ is a given deterministic function. Its well-posedness will be shown in Section 4.
Thus, there also exists a measurable map $\mathcal{G}^0:C([0,T];\mathbb{R})\rightarrow \mathcal{E}$ such that $\mathcal{G}^0(\int_{0}^{\cdot}h(s)ds)$ is the unique solution of (\ref{Skeleton}) for $h \in L^2([0,T];\mathbb{R})$.\par
    Our main theorem is stated as the follows.
    \begin{theorem}{}\label{main}
    	Under hypothesis (\hyperlink{H}{H}), the solution famility ${u^\varepsilon}$ of (\ref{eps}) satisfies  the LDP in $C([0,T];H)\cap L^2([0,T];V)$ with the rate function
    	$$I(f)=\inf\limits_{\{h\in L^2([0,T];\mathbb{R}):f=\mathcal{G}^0(\int_{0}^{\cdot}h(s)ds)\}}\Big\{\frac{1}{2}\int_{0}^{T}|h(s)|^2ds\Big\} $$
    	where $\mathcal{G}^0(\int_{0}^{\cdot}h(s)ds)$ is the unique solution of (\ref{Skeleton}) with $h \in L^2([0,T];\mathbb{R})$.
    \end{theorem}
     \noindent {\bf Proof}.
     The proof consists of two parts.\\
     \textbf{Part1:} Establish the  well-posedness of the skeleton equation (\ref{Skeleton}). This will be done in Section 4.
     \vskip 0.3cm
     \noindent \textbf{Part2:} Verification of the conditions in Theorem \ref{criterion}. This is done in Section 5, see Theorem \ref{a} and Theorem \ref{b} below.
     $\blacksquare$

\section{Skeleton equations}
\setcounter{equation}{0}
\quad In this section, we will prove the uniqueness and existence of solutions to equation (\ref{Skeleton}).  To do this, we will use the following estimates concerning to the logarithmic term, whose proof can be found in \cite{SZ}. In the following, we'll let $C$ denote an arbitrary constant, and  $C_T$ denote arbitrary constant depending on $T$, which could be different from line to line. Also, we will simply denote by  $\big\|\cdot\big\|_{L_T^2}$ the norm in $L^2([0,T];\mathbb{R})$.

\begin{lemma}\label{lemma 3.1}
For any $u,v\in V$, $\varepsilon>0$, and $\alpha\in(0,1)$, we have
\begin{align}\label{3.1}
  & \left(u\log|u| - v\log|v|, u-v\right) \nonumber\\
  \leq &\  \varepsilon \big\|u-v\big\|_{V}^2 + \left(1+\frac{d}{4}\log \frac{1}{\varepsilon} \right)  \big\|u-v\big\|^2
   + \big\|u-v \big\|^2\log\big\|u-v\big\|\nonumber\\  & + \frac{1}{2(1-\alpha)\mathrm{e}}\left(\big\|u\big\|^{2(1-\alpha)}+\big\|v\big\|^{2(1-\alpha)}\right)\big\|u-v\big\|^{2\alpha}.
\end{align}
\end{lemma}


\begin{lemma}\label{lemma 3.2}
  For any $u,v\in V$, $\varepsilon>0$, and $\alpha\in(0,1)$, we have
\begin{align}\label{20190629.1433.1}
  & \int_D |u(x)-v(x)|^2 \log_{+}\left(|u(x)|\vee |v(x)|\right) dx \nonumber\\
  \leq &\  \varepsilon \big\|u-v\big\|_{V}^2 + \left(\frac{d}{4}\log \frac{1}{\varepsilon} \right)  \big\|u-v\big\|^2
   + \big\|u-v \big\|^2\log\big\|u-v\big\|\nonumber\\  & + \frac{1}{2(1-\alpha)\mathrm{e}}\left(\big\|u\big\|^{2(1-\alpha)}+\big\|v\big\|^{2(1-\alpha)}\right)\big\|u-v\big\|^{2\alpha} \nonumber\\
   & + \frac{1}{2(1-\alpha)\mathrm{e}} \left(4m(D)\right)^{1-\alpha} \big\|u-v\big\|^{2\alpha},
\end{align}
where $m(D)$ is the Lebesgue measure of domain $D$.
\end{lemma}

\begin{theorem}\label{uniqueness}
Suppose hypothesis (\hyperlink{H}{H}) holds. Then the uniqueness of solutions holds for equation (\ref{Skeleton}).
\end{theorem}
\noindent {\bf Proof}.
 Fix $h\in L^2([0,T];\mathbb{R})$. Let $u^h, v^h\in L^2([0,T];V)\cap C([0,T];H)$ be two solutions of (\ref{Skeleton}).
Then
\begin{align*}
u^h(t)-v^h(t)=&\int_{0}^{t}\big(\Delta u^h(s)-\Delta v^h(s)\big)ds+\int_{0}^{t}\big(u^h(s)\log|u^h(s)|-v^h(s)\\&\log|v^h(s)|\big)ds+\int_{0}^{t}\big(\sigma(u^h(s))-\sigma(v^h(s))\big)h(s)ds.
\end{align*}
For $0<\delta<1$, we define the following time, with the convention that $\inf\{\emptyset\}=\infty$:
\begin{align*}
\tau^{\delta}:=&\inf\Big\{t\in(0,T]:\big\|u^h(t)-v^h(t)\big\|>\delta\Big\}.
\end{align*}
Set
\begin{align*}
M:=&\Big(\sup_{0\leq t\leq T}\big\|u^h(t)\big\|\Big)\vee \Big(\sup_{0\leq t\leq T}\big\|v^h(t)\big\| \Big)\vee \Big(\int_{0}^{T}\big\|u^h(s)\big\|_V^2ds \Big)\vee \Big(\int_{0}^{T}\big\|v^h(s)\big\|_V^2ds \Big).
\end{align*}
Let $Z_t:=u^h(t)-v^h(t)$. By the chain rule, we get
\begin{align}\label{220225.2040}
&\ \big\|Z_{t\wedge\tau^{\delta}}\big\|^2+2\int_{0}^{t\wedge\tau^{\delta}}\big\|Z_s\big\|_V^2ds\nonumber\\=&\ 2\int_{0}^{t\wedge \tau^{\delta}} \big(u^h(s)\log|u^h(s)|-v^h(s)\nonumber \log|v^h(s)|,u^h(s)-v^h(s)\big) ds\nonumber\\&+2\int_{0}^{t\wedge\tau^{\delta}} \big(\sigma(u^h(s))-\sigma(v^h(s)),u^h(s)-v^h(s)\big)h(s)ds\nonumber
	\\ \leq&\ \frac{1}{2} \int_{0}^{t\wedge\tau^{\delta}}\big\|Z_s\big\|_V^2ds+C\int_{0}^{t\wedge\tau^{\delta}}\big\|Z_s\big\|^2ds\nonumber\\&+\frac{1}{(1-\alpha)e}\int_{0}^{t\wedge\tau^{\delta}}
\Big(\big\|u^h(s)\big\|^{2(1-\alpha)}+\big\|v^h(s)\big\|^{2(1-\alpha)}\Big)\big\|Z_s\big\|^{2\alpha}ds\nonumber\\&+\int_{0}^{t\wedge\tau^{\delta}}\big\|Z_s\big\|^2h(s)^2ds
+\int_{0}^{t\wedge\tau^{\delta}}\big\|\sigma(u^h(s))-\sigma(v^h(s))\big\|^2ds\nonumber\\&+2\int_{0}^{t\wedge\tau^{\delta}}\big\|Z_s\big\|^2\log\big\|Z_s\big\|ds,
\end{align}
where the last inequality follows by taking $\varepsilon=\frac{1}{4}$ in Lemma \ref{lemma 3.1} and H\"{o}lder's inequality.
By (\hyperlink{H}{H}) and Lemma \ref{lemma 3.2} with $\varepsilon=\frac{1}{2L_2^2}$, we have
\begin{align}
	&\int_{0}^{t\wedge\tau^{\delta}}\big\|\sigma(u^h(s))-\sigma(v^h(s))\big\|^2ds\nonumber\\\leq&\ \frac{1}{2}\int_{0}^{t\wedge\tau^{\delta}}\big\|Z_s\big\|_V^2ds\nonumber+C\int_{0}^{t\wedge\tau^{\delta}}\big\|Z_s\big\|^2ds+C\int_{0}^{t\wedge\tau^{\delta}}\big\|Z_s\big\|^2\log\big\|Z_s\big\|ds\nonumber\\
&+\frac{C}{2(1-\alpha)e}\int_{0}^{t\wedge\tau^{\delta}}\left(\big\| u^h(s)\big\|^{2(1-\alpha)}+\big\| v^h(s)\big\|^{2(1-\alpha)}+(4m(D))^{1-\alpha}\right)\big\| Z_s\big\|^{2\alpha}ds.
\end{align}
Since $\delta<1$, we have
\begin{align}\label{220225.2041}
\int_{0}^{t\wedge\tau^{\delta}}\big\|Z_s\big\|^2\log\big\|Z_s\big\|ds\leq 0.
\end{align}
Now combining (\ref{220225.2040})-(\ref{220225.2041}) together and using the definition of $\tau^{\delta}$ we obtain
\begin{align*}
&\big\|Z_{t\wedge\tau^{\delta}}\big\|^2+\int_{0}^{t\wedge\tau^{\delta}}\big\|Z_s\big\|_V^2ds\\\leq&\int_{0}^{t\wedge\tau^{\delta}}C\big(1+|h(s)|^2\big)\big\|Z_s\big\|^2ds+\frac{C}{(1-\alpha)e}\int_{0}^{t\wedge\tau^{\delta}}
\big(C^{1-\alpha}+M^{2(1-\alpha)}\big)\big\|Z_s\big\|^{2\alpha}ds.
\end{align*}
Writing $Y_t:=\big\|Z_{t\wedge\tau^{\delta}}\big\|^2$, for some constant $C_1>0$ independent of $\alpha$, we have $$Y_t\leq\int_{0}^{t}C\big(1+|h(s)|^2\big)Y_sds+\frac{C^{1-\alpha}+M^{2-2\alpha}}{C_1(1-\alpha)}\int_{0}^{t}Y_s^{\alpha} ds.$$ Now by Lemma \ref{6.1} in the Appendix, and using H\"older's inequality, we have
\begin{align*}
	Y_t&\leq \left[\int_{0}^{t}\frac{C^{1-\alpha}+M^{2-2\alpha}}{ C_1}\exp\left((1-\alpha)\int_{s}^{t}C(1+|h(r)|^2)dr\right)ds\right]^{\frac{1}{1-\alpha}}\\&\leq\left(\frac{M^{2-2\alpha}}{C_1}+\frac{C^{1-\alpha}}{C_1}\right)^{\frac{1}{1-\alpha}}\left\{\int_{0}^{t}\exp\big(C_{T,\|h\|_{L^2_T}}\big)ds\right\}t^{\frac{\alpha}{1-\alpha}}\\&\leq\frac{1}{2}\left[(\frac{2t^\alpha}{C_1})^{\frac{1}{1-\alpha}}\times M^2+(\frac{2t^{\alpha}}{C_1})^{\frac{1}{1-\alpha}}\times C\right]\times C_{T,\|h\|_{L^2_T}}.
\end{align*}
Choosing $T^*$ small enough such that $T^*<\frac{C_1}{2}\wedge T$ and letting $\alpha\rightarrow 1$, we obtain $$Y_t=0,\  \forall\ t\in[0,T^*].$$
This implies that we must have $\tau^{\delta}>T^*$ and hence $Z_t=0,\  \forall\ t\in[0,T^*]$.
Since $T^*$ is independent of the initial value, starting from $T^*$ and repeating the same arguments, we can deduce that  $Z_t=0\ \text{on}\  [T^*,2T^*\wedge T]$. Continuing like this, we eventually get $Z_t=0\ \text{on}\  [0, T]$, proving the uniqueness.
$\blacksquare$
\\

Next, we establish the existence of the solution of the skeleton equation. To this end,
we first study the Galerkin approximating equations of (\ref{Skeleton}).

\vskip 0.3cm

Let $H_n$ denote the $n$-dimensional subspace of $H$ spanned by $\{e_1,\dots, e_n\}$. Let $P_n: V^*\rightarrow H_n$ be defined by
\begin{align}\label{4.1}
  P_n g := \sum_{i=1}^{n}\ \langle g,e_i\rangle e_i .
\end{align}
 Fix $h\in L^2([0,T];\mathbb{R})$. For any integer $n\geq 1$, we consider the following equation in the finite-dimensional space  $H_n$:
\begin{numcases}{}\label{galerkin}
  du_n(t)= \Delta u_n(t)dt+P_n [u_n(t)\log|u_n(t)|]dt+ P_n[\sigma(u_n(t))]h(t)dt,  \nonumber\\
 \ u_n(0)= P_n u_0 \label{4.2} .
\end{numcases}

Slightly modifying the proofs of Theorem 2.1 and Theorem 2.2 in \cite{FZ},  we arrive at the following results whose proof is omitted. 

\begin{theorem}
If hypothesis (\hyperlink{H}{H}) holds, then there exists a unique global solution $u_n$ to equation (\ref{4.2}).
\end{theorem}

The next step is to show that the Galerkin approximating solutions $\{u_n\}$  admits a limit $u$ which will be a solution of the skeleton equation (\ref{Skeleton}). To this end,  we will give some uniform estimates of the approximating solutions.

\begin{lemma}\label{galerkin moment}
	Under assumption (\hyperlink{H}{H}), the following estimate holds $$\sup\limits_n\left[\sup\limits_{t\in [0,T]}\big\|u_n(t)\big\|^2+\int_{0}^{T}\big\|u_n(s)\big\|_V^2ds\right]<\infty.$$
\end{lemma}
\noindent {\bf Proof}.
	From (\ref{galerkin}), it follows that
\begin{align}\label{220225.1401}
\big\|u_n(t)\big\|^2&=\big\|P_nu_0\big\|^2-2\int_{0}^{t}\big\|u_n(s)\big\|_V^2ds+2\int_{0}^{t}(u_n(s)\log|u_n(s)|,u_n(s))ds\nonumber\\&+2\int_{0}^{t}\big(\sigma(u_n(s)),u_n(s)\big)h(s)ds.
\end{align}
By (\ref{log sobolev}), we have
\begin{align}
\big(u_n(s)\log|u_n(s)|,u_n(s)\big)\leq \frac{1}{4}\big\|u_n(s)\big\|_V^2+(1+\frac{d}{4}\log4)\big\|u_n(s)\big\|^2+\big\|u_n(s)\big\|^2\log\big\|u_n(s)\big\|.
\end{align}
Using (\hyperlink{H}{H}) and (\ref{log sobolev modi}), we get
\begin{align}\label{220225.1402}
	\big|\big(\sigma(u_n(s)),u_n(s)\big)\big|\leq&\int_D\big|\sigma(u_n(s,x))\big|\big|u_n(s,x)\big|dx\nonumber\\\leq&\int_D\Big(L_3+L_4\big|u_n(s,x)\big|\big(\log_+\big|u_n(s,x)\big|\big)^{\frac{1}{2}}\Big)\big|u_n(s,x)\big|dx\nonumber\\\leq& \int_D \Big(C+C\big|u_n(s,x)\big|^2+L_4\big|u_n(s,x)\big|^2\log_+\big|u_n(s,x)\big|\Big)dx\nonumber\\\leq&\  C+C\big\|u_n(s)\big\|^2+L_4\Big(\theta(s)\big\|u_n(s)\big\|^2_V+\big(\frac{d}{4}\log\frac{1}{\theta(s)}\big)\big\|u_n(s)\big\|^2\nonumber\\&+\big\|u_n(s)\big\|^2\log\big\|u_n(s)\big\|^2\Big).
\end{align}
Here $\theta(\cdot):\mathbb{R}_+\rightarrow\mathbb{R}_+$ is any given  positive-valued function defined on the positive real axis.
Taking $\theta(s)=\frac{1}{4L_4(|h(s)|\vee1)}$, and combining (\ref{220225.1401})-(\ref{220225.1402}) together, we find
\begin{align*}
	&\ \big\|u_n(t)\big\|^2+\int_{0}^{t}\big\|u_n(s)\big\|_V^2ds\\\leq&\ \big\| u_0\big\| ^2+C_h+\int_{0}^{t}\big\| u_n(s)\big\|^2C\Big(1+|h(s)|+\log(|h(s)|\vee1)|h(s)|\Big)ds\\&\ +C\int_{0}^{t}\big\|u_n(s)\big\|^2\log\big\|u_n(s)\big\|^2\Big(1+|h(s)|\Big)ds.
\end{align*}
Then by Lemma \ref{6.2} in the Appendix, there exists some $C_{T,\|h\|_{L_T^2}}>0$ independent of $n$ such that,$$\big\|u_n(t)\big\|^2+\int_{0}^{t}\big\|u_n(s)\big\|_V^2ds
\leq(\big\|u_0\big\|^2+1)C_{T,\|h\|_{L_T^2}}.$$
This gives the desired estimate: $$\sup\limits_n\Big[\sup\limits_{t\in [0,T]}\big\|u_n(t)\big\|^2+\int_{0}^{T}\big\|u_n(s)\big\|_V^2ds\Big]<\infty. $$ $\blacksquare $\\


Let $W^{\beta, p}([0,T];V^*)$ be the space of measurable functionals $u(\cdot): [0,T]\rightarrow V^*$ with the finite norm defined by
\begin{align}\label{4.27}
	\big\| u\big\|^{p}_{{W^{\beta, p}([0,T];V^*)}}:= \int_0^{T} \big\| u(t)\big\|^p_{V^*} dt + \int_0^{T}\int_0^{T} \frac{\big\| u(t)-u(s)\big\|_{V^*}^p}{|t-s|^{1+\beta p}} dtds .
\end{align}
Next result shows  that $\{u_n\}_{n>0}$ is a bounded subset of $W^{\beta, p}([0,T];V^*)$ for some $\beta$ and $p$ .
\begin{lemma}\label{galerkin moment 2}
	Under hypothesis (\hyperlink{H}{H}), let $0<\beta<\frac{1}{2}$, then we have
	$$\sup\limits_n\Big
	\{\big\|u_n\big\|_{W^{\beta,2}([0,T];V^*)}\Big\}<\infty.$$
\end{lemma}
\noindent {\bf Proof}.
Since
\begin{align*}
	u_n(t)-u_n(s)=\int_{s}^{t}\Delta u_n(r)dr+\int_{s}^{t}P_n\Big[u_n(r)\log\big(u_n(r)\big)\Big]dr+\int_{s}^{t}P_n\Big[\sigma\big(u_n(r)\big)\Big]h(r)dr.		
\end{align*}
It follows that
\begin{align}\label{220225.1404}
\big\|u_n(t)-u_n(s)\big\|_{V^*}^2\leq&\ 3\times\Big[\Big(\int_{s}^{t}\big\|u_n(r)\big\|_Vdr\Big)^2+\Big(\int_{s}^{t}\big\|u_n(r)\log|u_n(r)|\big\|_{V^*}dr\ \Big)^2\nonumber\nonumber\\&+\Big(\int_{s}^{t}\big\|\sigma(u_n(r))\big\|_{V^*}h(r)dr\Big)^2\Big]\nonumber\\\leq&\  3\times(J_1+J_2+J_3).
\end{align}
By Lemma \ref{galerkin moment} and H\"older's inequality,
\begin{align}
 J_1\leq \left(\int_{s}^{t}\big\|u_n(r)\big\|_V^2dr\right)\big|t-s\big|.
\end{align}
By the Sobolev's embedding theorem and Riesz theorem,
$$L^{q^*}\hookrightarrow V^*,\ \forall\ q^*\geq\frac{2d}{d+2}.$$
Thus we can choose $\varepsilon>0,\ 2>q^*\geq\frac{2d}{d+2}$, such that when $d\geq3$,
\begin{numcases}{}
\big|x\log|x|\big|\leq C_\varepsilon(1+|x|^{1+\varepsilon}) \nonumber, \\
q^*(1+\varepsilon)\leq 2, \nonumber
\end{numcases}
for some $C_\varepsilon>0$.
Therefore,
\begin{align}
\big\|u_n(r)\log|u_n(r)|\big\|_{V^*}&\leq \left(\int_D\big(|u_n(r,x)|\log|u_n(r,x)|\big)^{q^*}dx\right)^{\frac{1}{q^*}}\nonumber\\&\leq  C+C\left(\int_D|u_n(r,x)|^{q^*(1+\varepsilon)}dx\right)^{\frac{1}{q^*}}\leq C+C\big\|u_n(r)\big\|^{1+\varepsilon}.\nonumber
\end{align}
Hence
\begin{align}
J_2\leq \Big(1+\sup\limits_n\sup\limits_{r\in[0,T]}\big\|u_n(r)\big\|^{1+\varepsilon}\Big)|t-s|^2.
\end{align}
Similarly, using the growth condition (\ref{H2}) it follows that
\begin{align}
	&\sup\limits_n\sup\limits_{r\in[0,T]}\big\|\sigma(u_n(r))\big\|_{V^*}<\infty\nonumber.
\end{align}
\vskip -0.3cm
\noindent
As a result,
\begin{align}\label{220225.1405}
J_3\leq C\Big[\int_{s}^{t}h(r)^2dr\Big]\big|t-s\big|.
\end{align}
Combining (\ref{220225.1404})-(\ref{220225.1405}) together and by Lemma \ref{galerkin moment}, for $\beta<\frac{1}{2}$, there exists a constant $C_{\|h\|_{L^2_T}}$, independent of $n$, such that
$$\big\|u_n\big\|_{W^{\beta,2}}^2=\int_{0}^{T}\big\|u_n(t)\big\|_{V^*}^2dt+\int_{0}^{T}\int_{0}^{T}\frac{\big\|u(t)-u(s)\big\|^2_{V^*}}{|t-s|^{1+2\beta}}dtds\leq C_{\|h\|_{L^2_T}}<\infty.$$
In the case when $d$=1 or 2, we can choose an arbitrary $q^*\in (1,2)$ and the same result follows.
$\blacksquare$\\

\vskip 0.3cm

Now we can establish the precompactness of the approximating solutions.
\begin{lemma}\label{compact}
	$\{u_n\}$ is precompact in $L^2([0,T];H)\cap C([0,T];V^*)$.
\end{lemma}
\noindent {\bf Proof}.
	From the proof of the Lemma \ref{galerkin moment 2} and the compact embedding $H\hookrightarrow V^*$, we conclude by the Arzela-Ascoli's theorem that  $\{u_n\}$ is precompact in $C([0,T];V^*)$.
On the other hand, it is known that $L^2([0,T];V)\cap W^{\beta,p}([0,T];V^*)$ is compactly embedded into $L^2([0,T];H)$ (see Theorem 4.5 in \cite{FG}). Lemma \ref{galerkin moment} and Lemma \ref{galerkin moment 2} together yield that  $\{u_n\}$ is bounded in $L^2([0,T];V)\cap W^{\beta,p}([0,T];V^*)$. So $\{u_n\}$ is also precompact in $L^2([0,T];H)$. $\blacksquare$

\vskip 0.3cm

Here is the main result of this section.
\begin{theorem}{}\label{existence}
	Suppose (\hyperlink{H}{H}) holds. Then for every initial value $u_0\in H$, there exists a unique solution to (\ref{Skeleton}).
\end{theorem}
\noindent {\bf Proof}.
By Lemma \ref{compact}, there exists a subsequence (still denoted by $\{u_n\}$) such that $u_n\rightarrow u$ in $L^2([0,T];H)\cap C([0,T];V^*)$.	Since $\big\|\cdot\big\|$ is lower-semicontinuous on $V^*$, we have for any $t\in[0,T]$,
\begin{align*}
	\big\|u(t)\big\|^2 \leq \liminf\limits_{n\rightarrow \infty}\big\|u_n(t)\big\|^2 \leq \sup_{n}\sup_{t\in [0,T]}\big\|u_n(t)\big\|^2 \leq C.
\end{align*}
Moreover, since  $u_n$ converges weakly (up to a subsequence) to $u$ in $L^2([0,T];V)$, it follows from Lemma \ref{galerkin moment} that
$$\int_{0}^{T}\big\|u(s)\big\|^2_Vds<\infty.$$
Furthermore, by Lemma \ref{galerkin moment} and the growth condition (\ref{H2}), we deduce the following results:
\begin{align*}
&(i)\ u_n(s,x)\rightarrow u(s,x) \quad   \text{a.e.}\ (s,x)\in[0,T]\times D,\\&
(ii)\ u_n\rightarrow u\ \text{weakly in }  L^2([0,T];V),\\&
(iii)\ \Delta u_n\rightarrow \Delta u\ \text{weakly in }  L^2([0,T];V^*),\\&
(iv)\ P_n\big[u_n\log|u_n|\big]\rightarrow u\log|u|\ \text{in}\ L^r([0,T];V^*)\ \text{for any } 1<r<2,\\&
(v)\ P_n\big[\sigma(u_n(\cdot))\big]\rightarrow \sigma(u(\cdot))\  \text{in} \ L^2([0,T];V^*).
\end{align*}
Now we let $n\rightarrow \infty$ in the equation satisfied by $u_n$ to see that $u$ satisfies
$$u(t)=u_0+\int_{0}^{t}\Delta u(s)ds+\int_{0}^{t}u(s)\log|u(s)|ds+\int_{0}^{t}\sigma(u(s))h(s)ds,$$
namely, $u$ is a solution of the skeleton equation.
The continuity of $u$ as an $H$-valued process follows from the above equation and the Lions-Magenes lemma (see e.g. Lemma 3.1.2 in \cite{T}). Therefore we obtain the existence of solutions to equation (\ref{Skeleton}). The uniqueness was proved in Theorem \ref{uniqueness}. $\blacksquare$

\section{Large Deviation Principle}
\setcounter{equation}{0}

\subsection{Moment estimate}
\quad  In this section, we will establish moment estimates for the solutions of (\ref{Abstracte}) and (\ref{Skeleton}) in preparation for the verification of the sufficient conditions stated  in Theorem \ref{criterion}. Throughout Section 5, we suppose that hypothesis (\hyperlink{H}{H}) is satisfied and $u_0\in H$.\\ \indent
For any fixed $N>0$ and for any family $\{h_\varepsilon\}_{\varepsilon>0}\subseteq \mathcal{A}_N$, by the Girsanov theorem and Yamada-Watanabe theorem (see \cite{RSZ}), $X^{h_\varepsilon}:=\mathcal{G}^\varepsilon(W.+\frac{1}{\varepsilon}\int_{0}^{\cdot}h_\varepsilon(s)ds)$ is the unique solution of the following equation,
\begin{numcases}{}\label{epsilon}
	X^{h_\varepsilon}(t) = u_0 +\int_0^t \Delta X^{h_\varepsilon}(s)ds + \int_0^t X^{h_\varepsilon}(s)\log|X^{h_\varepsilon}(s)| ds + \varepsilon\int_0^t \sigma(X^{h_\varepsilon}(s))dW_s \nonumber\\\qquad \quad \quad +\int_0^t \sigma(X^{h_\varepsilon}(s))h_\varepsilon(s)ds, \\
	\label{Abstractee}  X^{h_\varepsilon}(0)=u_0\in H. \nonumber
\end{numcases}
Moreover, $X^{h_\varepsilon}(\cdot)$ belongs to $\mathcal{E}(=C([0,T];H)\cap L^2([0,T];V)).$
For any fixed $N>0$, we denote by $Y^{h_\varepsilon}$ the solution of (\ref{Skeleton}) (the skeleton equation) with $h$ replaced by  $h_\varepsilon\in\mathcal{A}_N$, i.e $Y^{h_\varepsilon}_{\cdot}=\mathcal{G}^0(\int_{0}^{\cdot}h_\varepsilon(s)ds)$ according to the definition of the mapping $\mathcal{G}^0$. By Lemma \ref{galerkin moment} there exists a constant $C_N$ such that
\begin{align}\label{220225.2016}
\sup\limits_{\varepsilon>0}\Big\{\big\|Y^{h_\varepsilon}_\cdot\big\|_{L^\infty([0,T];H)}+\int_{0}^{T}\big\|Y^{h_\varepsilon}_s\big\|_V^2ds\Big\}<C_N, \ P\text{-a.s.}.
\end{align}
To prove the LDP,  we need a uniform moment estimate for $\{X^{h_\varepsilon}\}_{0<\varepsilon<1}$.
To this end, we define a family of stopping times:
$$\tau_M^{\varepsilon}=\inf\Big\{t\in (0,T]:\big\|X^{h_\varepsilon}_t\big\|\geq M\Big\},$$
$$\sigma_L^\varepsilon=\inf\Big\{t\in(0,T]:\int_{0}^{t}\big\|X^{h_\varepsilon}_s\big\|_V^2ds\geq L\Big\},$$
$$\tau^\varepsilon_{M,L}=\tau^\varepsilon_M\wedge\sigma^\varepsilon_L\wedge T,$$
 with the convention that $\inf\{\emptyset\}=\infty$.
We have $\tau_M^{\varepsilon}<\infty,\ \sigma^\varepsilon_L<\infty$, $P$-a.s by Girsanov's transformation. Inspired by Lemma 6.1 in \cite{SZ}, we define an auxiliary function $\Phi:[0,\infty)\rightarrow[0,\infty)$
$$\Phi(z):=\exp\left(\int_{0}^{z}\frac{dx}{1+x+x\rho(x)}\right),$$where
$$\rho(x)=
\begin{cases}
	\log x& x\geq e,\\
	\frac{x}{e}& x\in[0,e).
\end{cases}$$
Then $$\Phi'(z)=\Phi(z)\times\frac{1}{1+z+z\rho(z)},\quad \Phi''(z)\leq 0.$$
We have the following estimate.
\begin{theorem}\label{Th-1}For any$\ T,N>0$ and any family $\{h_\varepsilon\}_{\varepsilon>0}\subseteq\mathcal{A}_N$, there exists a constant $C_{T,N}$ such that
	 $$\sup\limits_{\varepsilon\in(0,1)\\ }E\Big[\sup\limits_{t\in[0,T]}\Phi(\big\|X^{h_\varepsilon}_{t}\big\|^2)+\int_{0}^{T}\Phi'(\big\|X^{h_\varepsilon}_s\big\|^2)\big\|X^{h_\varepsilon}_s\big\|_V^2ds\Big]\leq \Phi(\big\|u_0\big\|^2)C_{T,N}. $$
\end{theorem}

\noindent {\bf Proof}.
For any $L, M>0$ , by (\ref{epsilon}) and Ito's formula, we get
\begin{align*}
	&\ \big\|X_{t\wedge\tau^{\varepsilon}_{M,L}}^{h_\varepsilon}\big\|^2+2\int_{0}^{t\wedge\tau^{\varepsilon}_{M,L}}\big\|X^{h_\varepsilon}_s\big\|_V^2ds\\=&\ \|u_0\|^2+2\int_{0}^{t\wedge\tau^{\varepsilon}_{M,L}}\big(X^{h_\varepsilon}_s\log|X^{h_\varepsilon}_s|,X^{h_\varepsilon}_s\big)ds+2\varepsilon\int_{0}^{t\wedge\tau^{\varepsilon}_{M,L}}\big(\sigma(X^{h_\varepsilon}_s),X^{h_\varepsilon}_s\big)dW_s\\&+\varepsilon^2\int_{0}^{t\wedge\tau^{\varepsilon}_{M,L}}\big\|\sigma(X^{h_\varepsilon}_s)\big\|^2ds+2\int_{0}^{t\wedge\tau^{\varepsilon}_{M,L}}h_\varepsilon(s)\big(\sigma(X^{h_\varepsilon}_s),X^{h_\varepsilon}_s\big)ds.
\end{align*}
Applying again Ito's formula to the real-valued process $\big\|X_{t\wedge\tau^{\varepsilon}_{M,L}}^{h_\varepsilon}\big\|^2$, we obtain
\begin{align*}
  &\ \Phi(\big\|X_{t\wedge\tau^{\varepsilon}_{M,L}}^{h_\varepsilon}\big\|^2)+2\int_{0}^{t\wedge\tau^{\varepsilon}_{M,L}}\Phi'(\big\|X_{s}^{h_\varepsilon}\big\|^2)\big\|X^{h_\varepsilon}_s\big\|_V^2ds\\=&\  \Phi(\big\|u_0\big\|^2)+2\varepsilon\int_{0}^{t\wedge\tau^{\varepsilon}_{M,L}}\Phi'(\big\|X^{h_\varepsilon}_s\big\|^2)\big(\sigma(X^{h_\varepsilon}_s),X^{h_\varepsilon}_s\big)dW_s\\&+2\int_{0}^{t\wedge\tau^{\varepsilon}_{M,L}}\big(X^{h_\varepsilon}_s\log|X^{h_\varepsilon}_s|,X^{h_\varepsilon}_s\big)\Phi'\big(\big\|X^{h_\varepsilon}_s\big\|^2\big)ds\\&+2\varepsilon^2\int_{0}^{t\wedge\tau^{\varepsilon}_{M,L}}\Phi''\big(\big\|X^{h_\varepsilon}_s\big\|^2\big)\big(\sigma(X^{h_\varepsilon}_s),X^{h_\varepsilon}_s\big)^2ds\\&+\varepsilon^2\int_{0}^{t\wedge\tau^{\varepsilon}_{M,L}}\big\|\sigma(X^{h_\varepsilon}_s)\big\|^2\Phi'(\big\|X^{h_\varepsilon}_s\big\|^2)ds\\&+2\int_{0}^{t\wedge\tau^{\varepsilon}_{M,L}}h(s)\Phi'\big(\big\|X^{h_\varepsilon}_s\big\|^2\big)\big(\sigma(X^{h_\varepsilon}_s),X^{h_\varepsilon}_s\big)ds\\=&\text{ I+II+$\cdots$+VI}.
\end{align*}
It follows from Lemma \ref{lemma 3.1} and Lemma \ref{lemma 3.2} that
\begin{align}\label{220225.2035}
\text{III} \leq 2\int_{0}^{t\wedge\tau^{\varepsilon}_{M,L}}\Phi'\big(\big\|X^{h_\varepsilon}_s\big\|^2\big)\big(\frac{1}{8}\big\|X^{h_\varepsilon}_s\big\|^2_V+\frac{d\log8}{4}\big\|X^{h_\varepsilon}_s\big\|^2+\big\|X^{h_\varepsilon}_s\big\|^2\log\big\|X^{h_\varepsilon}_s\big\|\big)ds.
\end{align}
By the similar method as for (\ref{220225.1402}), there exists a positive constant $C_2$ such that
\begin{align}\label{220225.2034}
	\big\|\sigma\big(X^{h_\varepsilon}_s\big)\big\|^2&\leq C\Big(\int_{D}\big|X^{h_\varepsilon}_s(x)\big|^2\log_+\big|X^{h_\varepsilon}_s(x)\big|dx+1\Big)\nonumber\\
	&\leq C_2\left(\theta\big\|X^{h_\varepsilon}_s\big\|_V^2+\big(\log\frac{1}{\theta}+1\big)\big\|X^{h_\varepsilon}_s\big\|^2+\big\|X^{h_\varepsilon}_s\big\|^2\log\big\|X^{h_\varepsilon}_s\big\|+1\right),
\end{align}
where $\theta=\frac{1}{4C_2}$. Hence we find that
\begin{align}\label{220225.2036}
\text{V}\leq\int_{0}^{t\wedge\tau^{\varepsilon}_{M,L}}\varepsilon^2\Phi'(\big\|X^{h_\varepsilon}_s\big\|^2)\left(C+C\big\|X^{h_\varepsilon}_s\big\|^2\big(\log_+\big\|X^{h_\varepsilon}_s\big\|^2+1\big)+\frac{\big\|X^{h_\varepsilon}_s\big\|_V^2}{4}\right)ds.
\end{align}
Similarly as (\ref{220225.2034}), there exists a positive constant $C_3$ such that,
\begin{align}\label{220227.1741}
		\big|\big(\sigma(X^{h_\varepsilon}_s),X^{h_\varepsilon}_s\big)\big|\leq C_3\Big(\theta(s)\big\|X^{h_\varepsilon}_s\big\|_V^2+\big(\log\frac{1}{\theta(s)}+1\big)\big\|X^{h_\varepsilon}_s\big\|^2+\big\|X^{h_\varepsilon}_s\big\|^2\log\big\|X^{h_\varepsilon}_s\big\|+1\Big),
\end{align}
where $\theta(s)=\frac{1}{4C_3(|h(s)|\vee1)}$ in (\ref{220227.1741}). Hence we see that
\begin{align}\label{220225.2037}
\text{VI}\leq&\
 \frac{1}{2}\int_{0}^{t\wedge\tau^{\varepsilon}_{M,L}}\Phi'\big(\big\|X^{h_\varepsilon}_s\big\|^2\big)\big\|X^{h_\varepsilon}_s\big\|_V^2ds+C\int_{0}^{t\wedge\tau^{\varepsilon}_{M,L}}\Phi'\big(\big\|X^{h_\varepsilon}_s\big\|^2\big)\Big[\big(1+\big|h_\varepsilon(s)\big|^2\big)\big\|X^{h_\varepsilon}_s\big\|^2\nonumber\\&+|h_\varepsilon(s)|\big\|X^{h_\varepsilon}_s\big\|^2\log\big\|X^{h_\varepsilon}_s\big\|+\big|h_\varepsilon(s)\big|\Big]ds.
\end{align}
Since $\Phi''\leq0$, we have
\begin{align}\label{220225.2038}
\text{IV}\leq 0.
\end{align}
Combining (\ref{220225.2035}), (\ref{220225.2036}), (\ref{220225.2037}) and (\ref{220225.2038}) together, due to the fact that
\begin{align}\label{220225.1958}
1+x^2+x^2\log_+(x)\leq C_4(1+x^2+x^2\rho(x^2))\,
\end{align}
 for some $C_4>0$, we have
\begin{align*}
	&\ \Phi\big(\big\|X_{t\wedge\tau^{\varepsilon}_{M,L}}^{h_\varepsilon}\big\|^2\big)+\int_{0}^{t\wedge\tau^{\varepsilon}_{M,L}}\Phi'\big(\big\|X_{s}^{h_\varepsilon}\big\|^2\big)\big\|X^{h_\varepsilon}_s\big\|_V^2ds\\\leq&\ \Phi(\big\|u_0\big\|^2)+2\varepsilon\sup\limits_{r\in[0,t]}\Big|\int_{0}^{r\wedge\tau^{\varepsilon}_{M,L}}\Phi'(\big\|X^{h_\varepsilon}_s\big\|^2)(\sigma(X^{h_\varepsilon}_s),X^{h_\varepsilon}_s)dW_s\Big|\\&+\int_{0}^{t\wedge\tau^{\varepsilon}_{M,L}}\Phi'(\big\|X^{h_\varepsilon}_s\big\|^2)C\big(1+|h_\varepsilon(s)|^2\big)\big(1+\big\|X^{h_\varepsilon}_s\big\|^2\log_+\big\|X^{h_\varepsilon}_s\big\|+\big\|X^{h_\varepsilon}_s\big\|^2\big)ds\\\leq&\  \Phi\big(\big\|u_0\big\|^2\big)+2\varepsilon\sup\limits_{r\in[0,t]}\Big|\int_{0}^{r\wedge\tau^{\varepsilon}_{M,L}}\Phi'(\big\|X^{h_\varepsilon}_s\big\|^2)\big(\sigma(X^{h_\varepsilon}_s),X^{h_\varepsilon}_s\big)dW_s\Big|\\&+\int_{0}^{t\wedge\tau^{\varepsilon}_{M,L}}\Phi\big(\big\|X^{h_\varepsilon}_s\big\|^2\big)C\big(1+\big|h_\varepsilon(s)\big|^2\big)ds.
\end{align*}
By Gronwall's inequality, it follows that
\begin{align*}
&\sup\limits_{s\in[0,t]}\Phi\big(\big\|X^{h_\varepsilon}_{s\wedge\tau^{\varepsilon}_{M,L}}\big\|^2\big)+\int_{0}^{t\wedge\tau^{\varepsilon}_{M,L}}\Phi(\big\|X_{s}\big\|^2)\big\|X^{h_\varepsilon}_s\big\|_V^2ds\\\leq&\Big(\Phi\big(\big\|u_0\big\|^2\big)+2\sup\limits_{s\in[0,t]}\Big|\int_{0}^{s\wedge\tau^{\varepsilon}_{M,L}}\big(\sigma(X^{h_\varepsilon}_r),X^{h_\varepsilon}_r\big)\Phi'\big(\big\|X^{h_\varepsilon}_r\big\|^2)dW_r\Big|\Big)e^{C(t+N^2)}.
\end{align*}
Taking expectations on both sides of the above inequality, and then using  BDG's inequality, we obtain
\begin{align}\label{220225.2031}
&\ E\Big[\sup\limits_{s\in[0,t]}\Phi\big(\big\|X^{h_\varepsilon}_{s\wedge\tau^{\varepsilon}_{M,L}}\big\|^2\big)\Big]+E\Big[\int_{0}^{t\wedge\tau^{\varepsilon}_{M,L}}\Phi\big(\big\|X_{s\wedge\tau^{\varepsilon}_{M,L}}\big\|^2\big)\big\|X^{h_\varepsilon}_s\big\|_V^2ds\Big]\nonumber\\\leq&\left\{\Phi\big(\big\|u_0\big\|^2\big)+2E\Big[\Big(\int_{0}^{t\wedge\tau^{\varepsilon}_{M,L}}\big(\sigma(X^{h_\varepsilon}_s),X^{h_\varepsilon}_s\big)^2\Phi'\big(\big\|X^{h_\varepsilon}_s\big\|^2\big)^2ds\Big)^{\frac{1}{2}}\Big]\right\}e^{C(T+N^2)}.
\end{align}
Denote the constant $e^{C(T+N^2)}$ in the above inequality by $C_5$. We now give an estimate of the second term on the right hand side of the above inequality.
\begin{align}
	&\ E\Big[\Big(\int_{0}^{t\wedge\tau^{\varepsilon}_{M,L}}\big(\sigma(X^{h_\varepsilon}_s),X^{h_\varepsilon}_s\big)^2\Phi'\big(\big\|X^{h_\varepsilon}_s\big\|^2\big)^2ds\Big)^{\frac{1}{2}}\Big]\nonumber\\
\leq&\  E\Big[\Big(\int_{0}^{t\wedge\tau^{\varepsilon}_{M,L}}\big\|\sigma(X^{h_\varepsilon}_s)\big\|^2\big\|X^{h_\varepsilon}_s\big\|^2\Phi'\big(\big\|X^{h_\varepsilon}_s\big\|^2\big)^2ds\Big)^{\frac{1}{2}}\Big]\nonumber\\
\leq&\  E\Big[\sup_{s\in[0,t]}\Big(\big\|X^{h_\varepsilon}_{s\wedge\tau^{\varepsilon}_{M,L}}\big\|^2\Phi'\big(\big\|X^{h_\varepsilon}_{s\wedge\tau^{\varepsilon}_{M,L}}\big\|^2\big)\Big)^{\frac{1}{2}}\Big(\int_{0}^{t}
\big\|\sigma\big(X^{h_\varepsilon}_{s\wedge\tau^{\varepsilon}_{M,L}}\big)\big\|^2\Phi'\big(\big\|X^{h_\varepsilon}_{s\wedge\tau^{\varepsilon}_{M,L}}\big\|^2\big)ds\Big)^{\frac{1}{2}}\Big]\nonumber\\
\leq&\ \delta E\Big[\sup_{s\in[0,t]}\Big(\big\|X^{h_\varepsilon}_{s\wedge\tau^{\varepsilon}_{M,L}}\big\|^2\Phi'\big(\big\|X^{h_\varepsilon}_{s\wedge\tau^{\varepsilon}_{M,L}}\big\|^2\big)\Big)\Big]+\frac{1}{4\delta} E\Big[\int_{0}^{t}\big\|\sigma\big(X^{h_\varepsilon}_{s\wedge\tau^{\varepsilon}_{M,L}}\big)\big\|^2\Phi'\big(\big\|X^{h_\varepsilon}_{s\wedge\tau^{\varepsilon}_{M,L}}\big\|^2\big)ds\Big]\nonumber\\=&\ I_1+I_2,
\end{align}
 where $\delta$ remains to be chosen. Due to (\ref{220225.1958}), we have
\begin{align}
I_1\leq C_4\delta E\Big[\sup\limits_{s\in[0,t]}\Phi\big(\big\|X^{h_\varepsilon}_{s\wedge\tau^{\varepsilon}_{M,L}}\big\|^2\big)\Big].
\end{align}
For $I_2$, using (\ref{H2}), (\ref{log sobolev modi}) and Lemma \ref{lemma 3.2}, we get that
\begin{align}\label{220225.2032}
	&\ \big\|\sigma(X^{h_\varepsilon}_{s\wedge\tau^{\varepsilon}_{M,L}})\big\|^2\nonumber\\\leq&\int_{D}\Big(2L_3^2+2L_4^2\big|X^{h_\varepsilon}_{s\wedge\tau^{\varepsilon}_{M,L}}(x)\big|^2\log_+\big|X^{h_\varepsilon}_{s\wedge\tau^{\varepsilon}_{M,L}}(x)\big|\Big)dx\nonumber\\\leq&\  C+2L_4^2\varepsilon'\big\|X^{h_\varepsilon}_{s\wedge\tau^{\varepsilon}_{M,L}}\big\|_V^2+C_{\varepsilon'}\big\|X_{s\wedge\tau^{\varepsilon}_{M,L}}^{h_{\varepsilon}}\big\|^2+\big\|X^{h_\varepsilon}_{s\wedge\tau^{\varepsilon}_{M,L}}\big\|^2\log\big\|X^{h_\varepsilon}_{s\wedge\tau^{\varepsilon}_{M,L}}\big\|.
\end{align}
Here $\varepsilon'$ is also to be determined. Taking $\delta=\frac{1}{4C_4C_5},\  \varepsilon'=\frac{1}{8C_4C_5^2L_4^2}$ and combining (\ref{220225.2031})-(\ref{220225.2032}) together, we see that
\begin{align*}
&\ E\Big[\sup\limits_{s\in [0,t]}\Phi\big(\big\|X^{h_\varepsilon}_{s\wedge\tau^{\varepsilon}_{M,L}}\big\|^2\big)\Big]+E\Big[\int_{0}^{t\wedge\tau^{\varepsilon}_{M,L}}\Phi'\big(\big\|X^{h_\varepsilon}_s\big\|^2\big)\big\|X_s^{h_\varepsilon}\big\|_V^2ds\Big]\\ \leq&\  C_{T,N}\Phi(\big\|u_0\big\|^2)+C\int_{0}^{t}E\Big[\Phi\big(\big\|X^{h_\varepsilon}_{s\wedge\tau^{\varepsilon}_{M,L}}\big\|^2\big)\Big]ds.
\end{align*}
Using Gronwall's inequality,  we obtain
\begin{align*}
\sup\limits_{ \varepsilon\in (0,1)}{E\Big[\sup\limits_{s\in[0,T]}\Phi\big(\big\|X^{h_\varepsilon}_{s\wedge\tau^{\varepsilon}_{M,L}}\big\|^2\big)+\int_{0}^{T\wedge\tau^{\varepsilon}_{M,L}}\Phi'\big(\big\|X^{h_\epsilon}_{s}\big\|^2\big)\big\|X_s^{h_\varepsilon}\big\|_V^2ds\Big]}\leq C_{T,N}\Phi\big(\big\|u_0\big\|^2\big).
\end{align*}
Since $\tau^{\varepsilon}_{M,L}=\tau^\varepsilon_M\wedge\sigma^\varepsilon_L\wedge T \rightarrow T$ as $L, M\rightarrow \infty$, applying the Fatou lemma we get
$$\sup\limits_{\varepsilon\in(0,1) \\ }E\Big[\sup\limits_{t\in[0,T]}\Phi\big(\big\|X^{h_\varepsilon}_{t}\big\|^2\big)+\int_{0}^{T}\Phi'\big(\big\|X^{h_\varepsilon}_s\big\|^2\big)\big\|X^{h_\varepsilon}_s\big\|_V^2ds\Big]\leq \Phi\big(\big\|u_0\big\|^2\big)C_{T,N}. $$
$\blacksquare$

\vskip 0.5cm

\begin{corollary}\label{220225.1547}
	For any $T,N,M>0$ and any family $\{h_\varepsilon\}_{\varepsilon>0}\subseteq\mathcal{A}_N$, there exists a positive constant $C_{T,N,M}$, such that
	$$\sup\limits_{t\in[0,T]}\big\|X_{t\wedge \tau^\varepsilon_{M}}^{h_\varepsilon}\big\|\leq M, \quad  \sup\limits_{\varepsilon\in (0,1)  }E\Big[\int_{0}^{T\wedge\tau^\varepsilon_{M}}\big\|X_t^{h_\varepsilon}\big\|_V^2dt\Big]\leq C_{T,N,M}.$$
\end{corollary}
\noindent {\bf Proof}. From the definition of $\tau_M^\varepsilon$, it follows immediately  that $$\sup\limits_{t\in[0,T]}\big\|X_{t\wedge\tau^\varepsilon_{M}}^{h_\varepsilon}\big\|\leq M, \quad \forall\, \varepsilon\in (0,1).$$
	For $0\leq x\leq M$, $\Phi^{\prime}(x)\geq c_M$ for some $c_M>0$. It follows from Theorem \ref{Th-1} that
	\begin{align*}
	c_M \times E\Big[\int_{0}^{T\wedge\tau^\varepsilon_M}\big\|X_s^{h_\varepsilon}\big\|^2_Vds\Big]\leq E\Big[\int_{0}^{T\wedge\tau^\varepsilon_{M}}\Phi'(\big\|X_s^{h_\varepsilon}\big\|^2)\big\|X_s^{h_\varepsilon}\big\|_V^2ds\Big]\leq\Phi(\big\|u_0\big\|^2)C_{T,N}.
	\end{align*}
	Hence, $$\sup\limits_{\varepsilon\in (0,1)}E\Big[\int_{0}^{T\wedge\tau_{M}^\varepsilon}\big\|X_t^{h_\varepsilon}\big\|_V^2dt\Big]\leq C_{T,N,M}.$$
$\blacksquare$
\begin{corollary}\label{220225.2018}
		For any $T,N,M>0$ and any family $\{h_\varepsilon\}_{\varepsilon>0}\subseteq\mathcal{A}_N$, we have
	\begin{align}
	&\sup\limits_{\varepsilon\in (0,1)}P\big(\tau_M^\varepsilon\leq T\big)\rightarrow 0, \quad \text{as } M\rightarrow \infty,\\
\label{220225.2012}	&\sup\limits_{\varepsilon\in (0,1) }P\big(\sigma_L^\varepsilon\leq T\big)\rightarrow 0, \quad \text{as } L\rightarrow +\infty.
	\end{align}
\end{corollary}
\noindent {\bf Proof}. It follows from Theorem \ref{Th-1} that
 $$\sup\limits_{\varepsilon\in (0,1)}\Phi(M)P\big(\tau_M^\varepsilon\leq T\big)\leq\sup\limits_{\varepsilon\in (0,1)}E\Big[\Phi(\big\|X^{h_\varepsilon}_{\tau_M^\varepsilon\wedge T}\big\|^2)\Big]\leq\Phi\big(\big\|u_0\big\|^2\big)C_{T,N}.$$
 So,
 \begin{align}
 \sup\limits_{\varepsilon\in(0,1)}P(\tau_M^\varepsilon\leq T)\leq\frac{\Phi\big(\big\|u_0\big\|^2\big)C_{T,N}}{\Phi(M)}\rightarrow 0, \ as\ M\rightarrow \infty.
\end{align}
For any $M>0$, we have
 \begin{align}
 &\sup\limits_{\varepsilon\in (0,1) }P\big(\sigma_L^\varepsilon\leq T\big)\nonumber\\=&\sup\limits_{\varepsilon\in (0,1)}P\big(\int_{0}^{T}\big\|X_t^{h_\varepsilon}\big\|^2_Vdt> L\big)\nonumber\\\leq&\sup\limits_{\varepsilon\in (0,1)}\frac{1}{L}E\Big[\int_{0}^{T\wedge\tau_{M}^\varepsilon}\big\|X_t^{h_\varepsilon}\big\|_V^2dt\Big]+\sup\limits_{\varepsilon\in (0,1)}P\big(\tau_M^\varepsilon\leq T\big).
\end{align}
By Corollary \ref{220225.1547},  let first $L\rightarrow\infty$ and then $M\rightarrow\infty$ to get (\ref{220225.2012}).
$\blacksquare$
\subsection{Verification of condition (a) in Theorem \ref{criterion}}
\quad As a part of the proof of the main result Theorem 3.1, in this section, we will verify condition (a) in Theorem \ref{criterion}. We will first show that (a) holds on a small interval and then extend it to the whole interval $[0,T]$ by piecing  the small  intervals together.
Before the proof, we introduce the following notation for two functions $u(\cdot)$ and $v(\cdot)$.
$$\rho_{a,b}(u,v)^2=\sup\limits_{s\in[a,b]}\big\|u(s)-v(s)\big\|^2+\int_{a}^{b}\big\|u(s)-v(s)\big\|_V^2ds,$$
and write $\rho_a(u,v):=\rho_{0,a}(u,v)$. In the rest of the paper, we will denote the metric $\rho$ of space $\mathcal{E}$ by $\rho_T$.
\begin{theorem}{}\label{a}
For every $N<\infty$, any family $\{h_\varepsilon\}_{\varepsilon>0}\subseteq \mathcal{A}_N$ and any $\delta>0$,
\begin{equation}\label{condition-a}
\lim\limits_{\varepsilon\rightarrow 0}P\Big(\rho_T\big(Y^{h_\varepsilon},X^{h_\varepsilon}\big)>\delta\Big)=0
\end{equation}
where $X^{h_\varepsilon}:=\mathcal{G}^\varepsilon\big(W.+\frac{1}{\varepsilon}\int_{0}^{.}h_\varepsilon(s)ds\big)$, $Y^{h_\varepsilon}:=\mathcal{G}^0\big(\int_{0}^{.}h_\varepsilon(s)ds\big)$ and $\rho_T(\cdot,\cdot)$ stands for the metric of the space $\mathcal{E}=C([0,T];H)\cap L^2([0,T];V)$ defined above.
\end{theorem}
\noindent {\bf Proof}. We first prove (\ref{condition-a})  on a small interval.
\vskip 0.3cm
Let's define $$\tau^\varepsilon_{M,L,\delta}=\tau^\varepsilon_{M,L}\wedge\inf\big\{t\in(0,T]:\big\|X^{h_\varepsilon}_t-Y^{h_\varepsilon}_t\big\|>\delta\big\},$$ for any fixed $0<\delta<1$, with the convention that $\inf\{\emptyset\}=\infty$, where $\tau^\varepsilon_{M,L}=\tau^\varepsilon_{L}\wedge\tau^\varepsilon_{M}\wedge T$ was defined as at the beginning of Section 5. We will write $\tau^\varepsilon_{M,L,\delta}$ as $\tau^\varepsilon$ for simplicity.
By Ito's formula,
\begin{align}\label{220225.1628}
	&\ \big\|X_{t\wedge\tau^\varepsilon}^{h_\varepsilon}-Y_{t\wedge\tau^\varepsilon}^{h_\varepsilon}\big\|^2+2\int_{0}^{t\wedge\tau^\varepsilon}\big\|X_s^{h_\varepsilon}-Y_s^{h_\varepsilon}\big\|_V^2ds \nonumber\\
	=&\ 2\int_{0}^{t\wedge\tau^\varepsilon}\big(X_s^{h_\varepsilon}\log|X_s^{h_\varepsilon}|-Y_s^{h_\varepsilon}\log|Y_s^{h_\varepsilon}|,X_s^{h_\varepsilon}-Y_s^{h_\varepsilon}\big)ds \nonumber\\
	&+2\int_{0}^{t\wedge\tau^\varepsilon}\big(\sigma(X_s^{h_\varepsilon})-\sigma(Y_s^{h_\varepsilon}),X_s^{h_\varepsilon}-Y_s^{h_\varepsilon}\big)h_\varepsilon(s)ds \nonumber\\ &+2\varepsilon\int_{0}^{t\wedge\tau^\varepsilon}\big(\sigma(X_s^{h_\varepsilon}),X_s^{h_\varepsilon}-Y_s^{h_\varepsilon}\big)dW_s+\varepsilon^2\int_{0}^{t\wedge\tau^\varepsilon}\big\|\sigma(X_s^{h_\varepsilon})\big\|^2ds.
\end{align}
We will give estimates for each term. By Lemma \ref{3.1}, (\ref{220225.2016}) and the definition of $\tau^\varepsilon$, it follows that
\begin{align}\label{220225.1632}
	&\ 2\int_{0}^{t\wedge\tau^\varepsilon}\big(X_s^{h_\varepsilon}\log|X_s^{h_\varepsilon}|-Y_s^{h_\varepsilon}\log|Y_s^{h_\varepsilon}|,X_s^{h_\varepsilon}-Y_s^{h_\varepsilon}\big)ds\nonumber\\ \leq&\ \frac{1}{2}\int_{0}^{t\wedge\tau^\varepsilon}\big\|X_s^{h_\varepsilon}-Y_s^{h_\varepsilon}\big\|^2_Vds+C\int_{0}^{t\wedge\tau^\varepsilon}\big\|X_s^{h_\varepsilon}-Y_s^{h_\varepsilon}\big\|^2\big(1+\log\big\|X_s^{h_\varepsilon}-Y_s^{h_\varepsilon}\big\|\big)ds\nonumber\\&+\frac{2M^{2-2\alpha}}{2(1-\alpha)e}\int_{0}^{t\wedge\tau^\varepsilon}\big\|X_s^{h_\varepsilon}-Y_s^{\varepsilon}\big\|^{2\alpha}ds.
\end{align}
By the assumption (\hyperlink{H}{H}) and Lemma \ref{lemma 3.2}, there exists a constant $C_6>0$ such that for $s\leq \tau^\varepsilon$,
\begin{align}\label{5-1}
	&\ \big|\big(\sigma(X_s^{h_\varepsilon})-\sigma(Y_s^{h_\varepsilon}),X_s^{h_\varepsilon}-Y_s^{h_\varepsilon}\big)\big|\nonumber\\\leq&\
	C_6\int_{D}|X_s^{h_\varepsilon}(x)-Y_s^{h_\varepsilon}(x)|^2\big(1+\log_+(|X_s^{h_\varepsilon}|\vee |Y_s^{h_\varepsilon}|)\big)dx\nonumber\\\leq&\  C_6\Big(\big\|X_s^{h_\varepsilon}-Y_s^{h_\varepsilon}\big\|^2+\theta(s)\big\|X_s^{h_\varepsilon}-Y_s^{h_\varepsilon}\big\|^2_V+\log\frac{1}{\theta(s)}\big\|X_s^{h_\varepsilon}-Y_s^{h_\varepsilon}\big\|^2\nonumber\\
&+\big\|X_s^{h_\varepsilon}-Y_s^{h_\varepsilon}\big\|^2\log\big\|X_s^{h_\varepsilon}-Y_s^{h_\varepsilon}\big\|+\frac{2M^{2(1-\alpha)}+C^{1-\alpha}}{2(1-\alpha)e}\big\|X_s^{h_\varepsilon}-Y_s^{h_\varepsilon}\big\|^{2\alpha}\Big).
\end{align}
For $s\leq\tau^\varepsilon$, $\log\big\|X_s^{h_\varepsilon}-Y_s^{h_\varepsilon}\big\|\leq 0$. Taking $\theta(s)=\frac{1}{4C_6(|h_\varepsilon(s)|\vee1)}$, it follows from (\ref{5-1}) that there exists a positive constant $C_7$ such that
\begin{align}\label{220225.1631}
	&\ \Big|2\int_{0}^{t\wedge\tau^\varepsilon}\big(\sigma(X_s^{h_\varepsilon})-\sigma(Y_s^{h_\varepsilon}),X_s^{h_\varepsilon}-Y_s^{h_\varepsilon}\big)h_\varepsilon(s)ds\Big|\nonumber\\\leq&\ \frac{1}{2}\int_{0}^{t\wedge\tau^\varepsilon}\big\|X_s^{h_\varepsilon}-Y_s^{h_\varepsilon}\big\|_V^2ds+\int_{0}^{t\wedge\tau^\varepsilon}C(1+|h_\varepsilon(s)|^2) \big\|X_s^{h_\varepsilon}-Y_s^{h_\varepsilon}\big\|^2ds\nonumber\\&+\frac{M^{2(1-\alpha)}}{(1-\alpha)C_7}\int_{0}^{t\wedge\tau^\varepsilon}\big|h_\varepsilon(s)\big|\big\|X_s^{h_\varepsilon}-Y_s^{h_\varepsilon}\big\|^{2\alpha}ds.
\end{align}
Combining (\ref{220225.1628}), (\ref{220225.1632}) and (\ref{220225.1631}) together, we find
\begin{align*}
	&\sup\limits_{s\in[0,t]}\big\|X_{s\wedge\tau^\varepsilon}^{h_\varepsilon}-Y_{s\wedge\tau^\varepsilon}^{h_\varepsilon}\big\|^2+\int_{0}^{t\wedge\tau^\varepsilon}\big\|X_s^{h_\varepsilon}-Y_s^{h_\varepsilon}\big\|_V^2ds\\\leq&\int_{0}^{t}\big\|X_{s\wedge\tau^\varepsilon}^{h_\varepsilon}-Y_{s\wedge\tau^\varepsilon}^{h_\varepsilon}\big\|^2 C\big(1+|h_\varepsilon(s)|^2\big)ds+\frac{{M^{2(1-\alpha)}}}{(1-\alpha)C_7}\int_{0}^{t}\big\|X_{s\wedge\tau^\varepsilon}^{h_\varepsilon}-Y_{s\wedge\tau^\varepsilon}^{h_\varepsilon}\big\|^{2\alpha}  \big(1+|h_\varepsilon(s)|\big)ds\\&+C\varepsilon^2\int_{0}^{t\wedge\tau^\varepsilon}\big\|\sigma(X^{h_\varepsilon}_{s})\big\|^2ds+2\varepsilon\sup\limits_{r\in[0,t]}\Big|\int_{0}^{r\wedge\tau^\varepsilon}\big(\sigma(X_s^{h_\varepsilon}),X_s^{h_\varepsilon}-Y_s^{h_\varepsilon}\big)dW_s\Big|.
\end{align*}
By Lemma \ref{6.1} in the Appendix, the above implies that
\begin{align}\label{5-2}
&\sup\limits_{s\in[0,t]}\big\|X_{s\wedge\tau^\varepsilon}^{h_\varepsilon}-Y_{s\wedge\tau^\varepsilon}^{h_\varepsilon}\big\|^2+\int_{0}^{t\wedge\tau^\varepsilon}\big\|X_s^{h_\varepsilon}-Y_s^{h_\varepsilon}\big\|_V^2ds \nonumber\\
\leq & \Bigg\{\Bigg(C\varepsilon^2\int_{0}^{t}\big\|\sigma(X^{h_\varepsilon}_{s\wedge\tau^\varepsilon})\big\|^2ds+2\varepsilon\sup\limits_{r\in[0,t]}\Big|\int_{0}^{r\wedge\tau^\varepsilon}
\Big(\sigma(X_s^{h_\varepsilon}),X_s^{h_\varepsilon}-Y_s^{h_\varepsilon}\Big)dW_s\Big|\Bigg)^{1-\alpha} \nonumber\\ &\times \exp\Big((1-\alpha)\int_{0}^{t}C\big(1+|h_\varepsilon(s)|^2\big)ds\Big)+\int_{0}^{t}\frac{M^{2(1-\alpha)}}{C_7}\big(1+|h_\varepsilon(s)|\big) \nonumber\\
&\times \exp\Big((1-\alpha)C\int_{s}^{t}(1+|h_\varepsilon(r)|^2)dr\Big)ds\Bigg\}^{\frac{1}{1-\alpha}} \nonumber\\
\leq & \,2^{\frac{\alpha}{1-\alpha}}\Bigg\{\Bigg(C\varepsilon^2\int_{0}^{t}\big\|\sigma(X^{h_\varepsilon}_{s\wedge\tau^\varepsilon})\big\|^2ds+2\varepsilon\sup\limits_{r\in[0,t]}\Big|\int_{0}^{r\wedge\tau^\varepsilon}\Big(\sigma(X_s^{h_\varepsilon}),X_s^{h_\varepsilon}-Y_s^{h_\varepsilon}\Big)dW_s\Big|\Bigg)^{1-\alpha} \nonumber\\
&\times \exp\Big((1-\alpha)\int_{0}^{t}C\big(1+|h_\varepsilon(r)|^2\big)dr\Big)\Bigg\}^{\frac{1}{1-\alpha}}+2^{\frac{\alpha}{1-\alpha}}\Bigg\{\int_{0}^{t}\frac{M^{2(1-\alpha)}}{C_7}\big(1+|h_\varepsilon(s)|\big)\nonumber\\&\times \exp\Big((1-\alpha)C\int_{s}^{t}\Big(1+|h_\varepsilon(r)|^2\Big)dr\Big)ds\Bigg\}^{\frac{1}{1-\alpha}} \nonumber\\
=& \  2^{\frac{\alpha}{1-\alpha}}\Bigg\{\Big(C\varepsilon^2\int_{0}^{t}\big\|\sigma(X^{h_\varepsilon}_{s\wedge\tau^\varepsilon})\big\|^2ds+2\varepsilon\sup\limits_{r\in[0,t]}\big|\int_{0}^{r\wedge\tau^\varepsilon}
\big(\sigma(X_s^{h_\varepsilon}),X_s^{h_\varepsilon}-Y_s^{h_\varepsilon}\big)dW_s\big|\Big) \nonumber\\
&\times \exp\Big(\int_{0}^{t}C\big(1+|h_\varepsilon(s)|^2\big)ds\Big)\Bigg\} \nonumber\\
& + 2^{\frac{\alpha}{1-\alpha}}\Bigg\{\int_{0}^{t}\frac{M^{2(1-\alpha)}}{C_7}\big(1+|h_\varepsilon(s)|\big)
\exp\Big((1-\alpha)C\int_{s}^{t}\big(1+|h_\varepsilon(r)|^2\big)dr\Big)ds\Bigg\}^{\frac{1}{1-\alpha}}.
\end{align}
Without lose of generality, we assume $t\leq T$, then we have
\begin{align*}
	&\ \int_{0}^{t}\frac{M^{2(1-\alpha)}}{C_7}\Big(1+|h_\varepsilon(s)|\Big)\exp\Big((1-\alpha)C\int_{s}^{t}\Big(1+|h_\varepsilon(r)|^2\Big)dr\Big)ds\\\leq&\ \int_{0}^{t}\frac{M^{2(1-\alpha)}}{C_7}\Big(1+|h_\varepsilon(s)|\Big)\exp\big((1-\alpha)C_{N}\big)ds\\\leq&\ \Big(\int_{0}^{t}\big(1+|h_\varepsilon(s)|\big)^2ds\Big)^\frac{1}{2}t^{\frac{1}{2}}\frac{1}{C_7}\big(M^2\exp(C_{N})\big)^{1-\alpha}\\\leq&\  \big(t^{\frac{1}{2}}+N^\frac{1}{2}\big) t^{\frac{1}{2}}\times\frac{1}{C_7}(M^2C_{N})^{1-\alpha} \\
	\leq&\ \big(T^{\frac{1}{2}}+N^\frac{1}{2}\big) t^\frac{1}{2}\times\frac{1}{C_7}(M^2C_{N})^{1-\alpha}.
\end{align*}
Now, taking expectations on both sides of (\ref{5-2}), $\forall\, t\leq T$,
\begin{align*}
	&\ E\Big[\sup\limits_{r\in[0,t]}\big\|X_{r\wedge\tau^\varepsilon}^{h_\varepsilon}-Y_{r\wedge\tau^\varepsilon}^{h_\varepsilon}\big\|^2\Big]+E\Big[\int_{0}^{t\wedge\tau^\varepsilon}\big\|X_s^{h_\varepsilon}-Y_s^{h_\varepsilon}\big\|_V^2ds\Big]\\ \leq&\  C_{T,N}2^{\frac{\alpha}{1-\alpha}}E\Big[ C \varepsilon^2\int_{0}^{t}\big\|\sigma(X^{h_\varepsilon}_{s\wedge\tau^\varepsilon})\big\|^2ds+2\varepsilon\sup\limits_{r\in[0,t]}\Big|\int_{0}^{r\wedge\tau^\varepsilon}(\sigma(X_s^{h_\varepsilon}),X_s^{h_\varepsilon}-Y_s^{h_\varepsilon})dW_s\Big|\Big] \\ & + \Big(2^\alpha(T^\frac{1}{2}+N^\frac{1}{2})\Big)^{\frac{1}{1-\alpha}} t^\frac{1}{2(1-\alpha)}\times\frac{1}{C_7^{\frac{1}{1-\alpha}}}M^2C_{N}.
\end{align*}
Let now $\varepsilon\rightarrow 0$ to obtain that
\begin{align*}
	&\ \limsup\limits_{\varepsilon\rightarrow 0}E\Big[\big\|X_{t\wedge\tau^\varepsilon}^{h_\varepsilon}-Y_{t\wedge\tau^\varepsilon}^{h_\varepsilon}\big\|^2\Big]+E\Big[\int_{0}^{t\wedge\tau^\varepsilon}\big\|X_s^{h_\varepsilon}-Y_s^{h_\varepsilon}\big\|_V^2ds\Big]\\\leq&\ M^2C_N\Big[2^\alpha(T^\frac{1}{2}+N^\frac{1}{2})t^\frac{1}{2}\times\frac{1}{C_7}\Big]^{\frac{1}{1-\alpha}}.
\end{align*}
We take $T^*=\frac{C_7^2}{16(T+N)}\wedge T$ and let $\alpha\rightarrow 1$ to see that  for $t\leq T^*$,
\begin{align}\label{220225.1650}
\lim\limits_{\varepsilon\rightarrow 0}E\Big[\big\|X_{t\wedge\tau^\varepsilon}^{h_\varepsilon}-Y_{t\wedge\tau^\varepsilon}^{h_\varepsilon}\big\|^2\Big]+E\Big[\int_{0}^{t\wedge\tau^\varepsilon}\big\|X_s^{h_\varepsilon}-Y_s^{h_\varepsilon}\big\|_V^2ds\Big]=0.
\end{align}
By the definition of $\tau^\varepsilon=\tau^\varepsilon_{M,L,\delta}$ and Chebyshev's inequality,
\begin{align}
&P\Big(\sup_{t\in[0,T^*]}\big\|X_t^{h_\varepsilon}-Y_t^{h_\varepsilon}\big\|>\delta\Big)\\\leq&\sup\limits_{\varepsilon\in (0,1) }P\Big(\tau_L^\varepsilon\leq T^*\Big)+\sup\limits_{\varepsilon\in(0,1)}P\Big(\tau_M^\varepsilon\leq T^*\Big)\nonumber +\frac{1}{\delta^2}E\Big[\big\|X_{T^*\wedge\tau^\varepsilon}^{h_\varepsilon}-Y_{T^*\wedge\tau^\varepsilon}^{h_\varepsilon}\big\|^2\Big],
\end{align}
 which tends to $0$ as $\varepsilon\rightarrow 0$ by Corollary \ref{220225.2018} and (\ref{220225.1650}). Similarly, we also can see that
\begin{align}
	& P\Big(\int_{0}^{T^*}\big\|X_s^{h_\varepsilon}-Y_s^{h_\varepsilon}\big\|_V^2ds>\delta^2\Big) \longrightarrow 0, \quad\text{as } \varepsilon\rightarrow 0.
\end{align}
Hence
\begin{align}	
	&\ P\Big(\rho_{T^*}(X^{h_\varepsilon},Y^{h_\varepsilon})^2>2\delta^2\Big)\nonumber\\\leq&\  P\Big(\int_{0}^{T^*}\big\|X_s^{h_\varepsilon}-Y_s^{h_\varepsilon}\big\|_V^2ds>\delta^2\Big)+P\Big(\sup_{t\in[0,T^*]}\big\|X_t^{h_\epsilon}-Y_t^{h_\epsilon}\big\|>\delta\Big)\nonumber\\
	& \longrightarrow 0 , \quad \text{as } \varepsilon\rightarrow 0.
	\end{align}
Thus, we have verified condition (a) of Theorem \ref{criterion} on a small interval $[0,T^*]$.  Now consider the equations satisfied by $X^{h_\varepsilon}$ and $Y^{h_\varepsilon}$ on the interval $[T^*, T]$ with respectively the initial values $X_{T^*}^{h_\varepsilon}$ and $Y_{T^*}^{h_\varepsilon}$.  Since $\big\|X_{T^*}^{h_\varepsilon}-Y_{T^*}^{h_\varepsilon}\big\|\rightarrow 0$ in probability $P$ as $\varepsilon\rightarrow 0$, using the same arguments as above, we can show that
\begin{align}
	&P\Big(\rho_{T^*,2T^*\wedge T}(X^{h_\varepsilon},Y^{h_\varepsilon})^2>2\delta^2\Big)\rightarrow 0,
\end{align}
and similarly for $n\geq 2$,
\begin{align}
	&P\Big(\rho_{(n-1)T^*,nT^*\wedge T}(X^{h_\varepsilon},Y^{h_\varepsilon})^2>2\delta^2\Big)\rightarrow 0.
\end{align}
Since there exists some $n>0$ such that $nT^*\geq T$, we finally obtain that
$$P\Big(\rho_T(X^{h_\varepsilon},Y^{h_\varepsilon})^2>\delta\Big)\rightarrow 0,$$
as $\varepsilon\rightarrow0$. $\blacksquare$

\subsection{Verification of condition (b) in Theorem \ref{criterion}}
\quad In this section, we will show that condition (b) of Theorem \ref{criterion} holds. Recall  that
 \begin{numcases}{}\label{5-4}
    	Y^{h_\varepsilon}(t) = u_0 +\int_0^t \Delta Y^{h_\varepsilon}(s)ds + \int_0^t Y^{h_\varepsilon}(s)\log|Y^{h_\varepsilon}(s)| ds + \int_0^t {h_\varepsilon}(s)\sigma(Y^{h_\varepsilon}(s))ds , \nonumber\\
    	 Y^{h_\varepsilon}(0)=u_0\in H,
    \end{numcases}
 $Y^{h}$ denotes the solution of (\ref{5-4}) with ${h_\varepsilon}$ replaced by $h$. Recall also
 $$\mathcal{E}=C([0,T];H)\cap L^2([0,T];V).$$
\vskip 0.4cm
The following theorem shows that  condition (b) in Theorem \ref{criterion} is indeed satisfied.
\vskip 0.5cm
\begin{theorem}{}\label{b}
	For every $N<\infty$ and any sequence $\{h_\varepsilon\}_{\varepsilon>0}\subseteq S_N$, if $h_\varepsilon$ converges to some element $h$ weakly in $L^2([0,T];\mathbb{R})$ as $\varepsilon\rightarrow 0$, then we have $Y^{h_\varepsilon}\rightarrow Y^{h}\ in\ \mathcal{E}$.
\end{theorem}
\noindent {\bf Proof.}
	From the proof of  Lemma \ref{galerkin moment}, we can see that
\begin{align}\label{5-5}
\sup\limits_{\varepsilon}\Big[\sup\limits_{t\in [0,T]}\big\|Y^{h_\varepsilon}_t\big\|^2
+\int_{0}^{T}\big\|Y^{h_\varepsilon}_s\big\|_V^2ds\Big]<\infty.
\end{align}
Using the same arguments as in the proofs of Lemma \ref{galerkin moment 2} and Lemma \ref{compact}, we can prove that $\{Y^{h_\varepsilon}\}_{\varepsilon>0}$ is a precompact subset of $L^2([0,T];H)\cap C([0,T];V^*)$. Then there exists a sequence $\varepsilon_n\rightarrow 0$ and some $v\in L^2([0,T];H)\cap C([0,T];V^*)$  such that
\begin{align*}
	\sup\limits_{t\in[0,T]}\big\|Y_t^{h_{\varepsilon_n}}-v(t)\big\|_{V^*}+ \int_{0}^{T}\big\|Y_s^{h_{\varepsilon_n}}-v(s)\big\|^2ds \xrightarrow[n\rightarrow\infty]{} 0.
\end{align*}
Also, by the lower semi-continuity of the corresponding norms, it follows from (\ref{5-5}) that
\begin{align}\label{220225.1727}
\sup\limits_{t\in[0,T]}\big\|v(t)\big\|<\infty\ \ \text{and }
\int_{0}^{T}\big\|v(t)\big\|_V^2dt<\infty.
\end{align}
Now we'll show that $v=Y^h$. Using the Sobolev imbedding, up to a subsequence,  the following holds:
\begin{align*}
	&(i)\ Y^{h_{\varepsilon_n}}\rightarrow v\ \text{in}\ C([0,T];V^*)\cap L^2([0,T];H),\\
	&(ii)\ Y^{h_{\varepsilon_n}}\rightarrow v\ \text{weakly in}\   L^2([0,T];V),\\
	&(iii)\ \Delta Y^{h_{\varepsilon_n}}\rightarrow \Delta v\  \text{weakly in}\ L^2([0,T];V^*),\\
	&(iv)\ Y^{h_{\varepsilon_n}}\log|Y^{h_{\varepsilon_n}}|\rightarrow \ v\log|v|\ \text{in} \ L^r([0,T];V^*)\text{ for any $1<r<2$},\\
	&(v)\ \sigma(Y^{h_{\varepsilon_n}}_\cdot)\rightarrow\sigma(v(\cdot))\  \text{in}\ L^2([0,T];V^*).
\end{align*}
Passing to the limit in the equation satisfied by $Y^{h_{\varepsilon_n}}$ as $n\rightarrow\infty$, we see that $v$ is a solution of the equation  (\ref{5-4}) with $h_\varepsilon$ replaced by $h$ using the weak convergence of ${h_\varepsilon}$. Due to the uniqueness of the solution of the skeleton equation, we have  $v=Y^h$. Because the limit point $v$ is unique, this enables us to conclude that
\begin{equation}\label{5-6}
\sup\limits_{t\in[0,T]}\big\|Y_t^{h_{\varepsilon}}-Y_t^h\big\|_{V^*}+ \int_{0}^{T}\big\|Y_s^{h_{\varepsilon}}-Y_s^h\big\|^2ds \xrightarrow[\varepsilon\rightarrow 0]{} 0 .
\end{equation}
To complete the proof, it remains to show that  $\{Y^{h_\varepsilon}\}$ actually converges  to $Y^h$ under the stronger metric $\rho_T(\cdot, \cdot)$.
Apply the chain rule and Lemma \ref{lemma 3.1} to find that
\begin{align*}
	&\ \big\|Y_t^{h_{\varepsilon}}-Y_t^h\big\|^2+2\int_{0}^{t}
\big\|Y_s^{h_{\varepsilon}}-Y_s^{h}\big\|_V^2ds \\
=&\int_{0}^{t}(Y_s^{h_{\varepsilon}}\log|Y^{h_{\varepsilon}}_s|
-Y_s^{h}\log|Y^{h}_s|,Y_s^{h_\varepsilon}-Y_s^{h})ds
+\int_{0}^{t}\Big(\sigma(Y_s^{h_{\varepsilon}})h_{\varepsilon}(s)-\sigma(Y_s^{h})h(s),
Y^{h_{\varepsilon}}_s-Y^h_s\Big)ds \\
	\leq & \int_{0}^{t}\big\|Y_s^{h_{\varepsilon}}-Y_s^{h}\big\|_V^2ds +C\int_{0}^{t}\big\|Y_s^{h_{\varepsilon}}-Y_s^{h}\big\|^2ds +\int_{0}^{t}\big\|Y_s^{h_{\varepsilon}}-Y_s^{h}\big\|^2\log\big\|Y_s^{h_{\varepsilon}}-Y_s^{h}\big\|ds \\ & +C_\alpha\int_{0}^{t}\Big(\big\|Y_s^{h_{\varepsilon}}\big\|^{2(1-\alpha)}+\big\|Y_s^{h}\big\|^{2(1-\alpha)}\Big)\big\|Y_s^{h_{\varepsilon}}-Y_s^{h}\big\|^{2\alpha}ds  \\ & +\int_{0}^{t}\Big(\sigma\big(Y_s^{h_{\varepsilon}}\big)h_{\varepsilon}(s)-\sigma\big(Y_s^{h}\big)h(s),Y^{h_{\varepsilon}}_s-Y^h_s\Big)ds.
\end{align*}
Therefore,
\begin{align*}
&\ \big\|Y_t^{h_{\varepsilon}}-Y_t^h\big\|^2+\int_{0}^{t}
\big\|Y_s^{h_{\varepsilon}}-Y_s^{h}\big\|_V^2ds \\
\leq &\ C\int_{0}^{t}\big\|Y_s^{h_{\varepsilon}}-Y_s^{h}\big\|^2ds + \int_{0}^{t}\big\|Y_s^{h_{\varepsilon}}-Y_s^{h}
\big\|^2\log_+\big\|Y_s^{h_{\varepsilon}}-Y_s^{h}\big\|ds \\
& +C_\alpha\int_{0}^{t}\Big(\big\|Y_s^{h_{\varepsilon}}\big\|^{2(1-\alpha)}
+\big\|Y_s^{h}\big\|^{2(1-\alpha)}\Big)\big\|Y_s^{h_{\varepsilon}}-Y_s^{h}
\big\|^{2\alpha}ds \\
& +\int_{0}^{t}\Big(\sigma(Y_s^{h_{\varepsilon}})h_{\varepsilon}(s)
-\sigma(Y_s^{h})h(s),Y^{h_{\varepsilon}}_s-Y^h_s\Big)ds\\
=&\ \text{I+II+III+IV}.
\end{align*}
We will finish the proof of the Theorem if  we show that the right hand side converges to 0, uniformly over $t\in[0,T]$. Since $Y^{h_{\varepsilon}}\rightarrow Y^h\ \text{in}\ C([0,T];V^*)\cap L^2([0,T];H)$,  taking into account (\ref{5-5}) we see that $\text{I+II+III}\rightarrow 0$ uniformly over $t\in[0,T]$. Write
\begin{align*}
	\text{IV}=&\int_{0}^{t}\Big(\sigma(Y_s^{h_{\varepsilon}})h_{\varepsilon}(s)
-\sigma(Y_s^{h})h(s),Y^{h_{\varepsilon}}_s-Y^h_s\Big)ds\\
=&\int_{0}^{t}\big(h_{\varepsilon}(s)-h(s)\big)\Big(\sigma(Y_s^{h}),
Y_s^{h_{\varepsilon}}-Y_s^{h}\Big)ds\\&+\int_{0}^{t}h_{\varepsilon}(s)
\Big(\sigma(Y_s^{h_{\varepsilon}})-\sigma(Y_s^{h}),Y_s^{h_{\varepsilon}}
-Y_s^{h}\Big)ds
	\\=&\ \text{IV}_1+\text{IV}_2.
\end{align*}
The estimate for $\text{IV}_1$ is as follows. For any $M>0$, let
\begin{align*}
 A_M:=\{s\in[0,T]:\big\|\sigma(Y^h_s)\big\|\leq M\}, \quad A_M^c=[0,T]\setminus{A_M}.
\end{align*}
 Obviously, for each $M>0$,
$$\int_{0}^{t}\big(\sigma(Y^h(s))1_{A_M},Y_s^{h_{\varepsilon}}-Y_s^{h}\big)^2ds\leq M^2\int_{0}^{T}\big\|Y_s^{h_{\varepsilon}}-Y_s^{h}\big\|^2ds\rightarrow 0, \ \   \text{as } \varepsilon\rightarrow 0,$$
which yields that
 $$\int_{0}^{t}\big(h_{\varepsilon}(s)-h(s)\big)\big(\sigma(Y_s^{h})
 1_{A_M},Y_s^{h_{\varepsilon}}-Y_s^{h}\big)ds\rightarrow0,$$
 uniformly over $t\in[0,T]$ as $\varepsilon\rightarrow 0$.  Hence to show $\text{IV}_1\rightarrow 0$, it suffices to prove that
 \begin{align}\label{220225.1744}
 \sup_\varepsilon\Big|\int_{0}^{t}\big(h_{\varepsilon}(s)-h(s)\big)
 \big(\sigma(Y^h(s))1_{A_M^c},Y_s^{h_{\varepsilon}}-Y_s^{h}\big)ds
 \Big|\rightarrow 0, \quad \text{as } M\rightarrow\infty,
 \end{align}
uniformly over $t\in[0,T]$.
By H\"older's inequality,
 \begin{align}\label{220225.1747}
 	&\sup_\varepsilon\Big|\int_{0}^{t}\big(h_{\varepsilon}(s)-h(s)\big)
 \big(\sigma(Y^h(s))1_{A_M^c},Y_s^{h_{\varepsilon}}-Y_s^{h}\big)ds
 \Big| \nonumber\\
 \leq&\sup_\varepsilon\int_{0}^{T}\big|h_{\varepsilon}(s)-h(s)\big|
 \big\|\sigma(Y_s^h)\big\|_{L^r(D)}\big\|Y_s^{h_{\varepsilon}}
 -Y_s^h\big\|_{L^{r'}(D)}1_{A_M^c}ds \nonumber\\
 \leq&\sup_\varepsilon\sup_{s\leq T}\big\|\sigma(Y_s^h)\big\|_{L^r(D)}\Big(\int_{0}^{T}
 \big|h_{\varepsilon}(s)-h(s)\big|^2ds\Big)^\frac{1}{2}
 \Big(\int_{0}^{T}\big\|Y_s^{h_{\varepsilon}}-Y_s^h\big\|^\alpha_{L^{r'}(D)}ds
 \Big)^\frac{1}{\alpha}m(A_M^c)^{\frac{1}{\alpha'}},
 \end{align}
where $\alpha$,$\alpha'$,$r$,$r'$ are positive real numbers such that
$$\frac{1}{\alpha}+\frac{1}{\alpha'}+\frac{1}{2}=1, \ \  \frac{1}{r}+\frac{1}{r'}=1, \ \text{and }\ 1<r<2.$$
By (\ref{220225.1727}) and (\ref{log sobolev modi}), there exists a positive constant $C_N$ such that
\begin{align}
	 \sup_{s\leq T}\big\|\sigma(Y_s^h)\big\|_{L^{r}(D)}\leq C_N, \quad \int_{0}^{T}\big\|\sigma(Y_s^{h})\big\|^2ds\leq C_N.
\end{align}
Hence
\begin{align}\label{220225.1748}
m(A_M^c)\leq\frac{1}{M^2}\int_{0}^{T}\|\sigma(Y_s^h)\|^2ds\leq\frac{C_N}{M^2}.
\end{align}
Combining (\ref{220225.1747})-(\ref{220225.1748}) together, to prove (\ref{220225.1744}), it remains to show that \begin{align}
\sup\limits_{\varepsilon}\int_{0}^{T}\big\|Y^{h_{\varepsilon}}_s-Y^h_s\big\|_{L^{r'}(D)}^\alpha ds<\infty, \quad\text{for some } \alpha, r' >2.
\end{align}
In fact, applying the interpolation inequality, Sobolev's embedding inequality, (\ref{5-5}) and (\ref{220225.1727}), we can find positive constants $\alpha,\ r',\ p\in (2,\infty)$ and $\delta\in(0,1)$ with
$$
\begin{cases}
	\frac{1}{r'}=\frac{\delta}{2}+\frac{1-\delta}{p}, \\
	\frac{1}{p}=\frac{1}{2}-\frac{1}{d}, \\
	\alpha(1-\delta)=2, \\
\end{cases}
$$
such that
\begin{align}\label{220225.1807}
	\sup\limits_\varepsilon\int_{0}^{T}\big\|Y^{h_{\varepsilon}}_s-Y^h_s
\big\|_{L^{r'}(D)}^\alpha ds &\leq\sup\limits_{\varepsilon}\int_{0}^{T}\big\|Y^{h_{\varepsilon}}_s
-Y^h_s\big\|^{\alpha\delta}\big\|Y^{h_{\varepsilon}}_s-Y^h_s
\big\|_{L^{p}(D)}^2\nonumber\\&\leq C_N\sup\limits_\varepsilon\int_{0}^{T}\big\|Y^{h_{\varepsilon}}_s-Y^h_s\big\|_{V}^2 ds<\infty.
\end{align}
Such constants $\alpha,r',\delta,p$ exist when $d\geq3$. For $d=1$ or $2$, we simply take $p=3$, and
$$
\begin{cases}
	\frac{1}{r'}=\frac{\delta}{2}+\frac{1-\delta}{3}, \\
	\alpha(1-\delta)=2. \\
\end{cases}
$$
Therefore (\ref{220225.1744}) holds. As a result, $\text{IV}_1\rightarrow 0$ uniformly over $t\in [0,T]$ as $\varepsilon\rightarrow 0$.

Now we turn to the convergence of $\text{IV}_2$, by (\ref{5-5}) and (\ref{220225.1727}) it suffices to show that
\begin{equation}\label{5-7}
\int_{0}^{T}\big\|\sigma(Y_s^{h_{\varepsilon}})-\sigma(Y_s^{h})\big\|^2ds
=\int_{0}^{T}\int_{D}\big|\sigma(Y_s^{h_{\varepsilon}}(x))
-\sigma(Y_s^{h}(x))\big|^2dxds\rightarrow 0.
\end{equation}
Since $Y_s^{h_{\varepsilon}}(x)
\rightarrow Y_s^{h}(x)$ for a.e. $(s,x)$, it suffices to show that for some $\alpha>2$,
$$
\sup\limits_\varepsilon\int_{0}^{T}\int_{D}\big|\sigma(Y_s^{h_{\varepsilon}}(x))
-\sigma(Y_s^{h}(x))\big|^\alpha dxds<\infty.
$$
By the growth condition on $\sigma$, it suffices to prove that for some $\beta>2$
\begin{align}\label{220225.1810}
\sup\limits_\varepsilon\int_{0}^{T}\int_{D}\big|Y_s^{h_{\varepsilon}}(x)\big|^\beta dxds<\infty.
\end{align}
Now taking $\beta=\frac{2d+4}{d}$, using the interpolation inequality, Sobolev's embedding inequality and (\ref{5-5}), as the proof of (\ref{220225.1807}), we can see that (\ref{220225.1810}) holds.
$\blacksquare$

\vskip 0.5cm
\section{Appendix}
\quad In this section, introduce two nonlinear-type Gronwall's inequality used in the paper.
\begin{lemma}\label{6.1}
	Let a,b,Y be nonnegative functions on $\mathbb{R}_+$, and there are constants $c\geq 0$, $0\leq\alpha<1$ such that
	 $$Y(t)\leq c+\int_{t_0}^{t}\Big(a(s)Y(s)+b(s)Y(s)^\alpha\Big)ds,\ \forall\ t\geq t_0\geq0.$$
	Then for any $t\geq t_0$,
\begin{align*}
Y(t)&\leq\Bigg{\{}c^{1-\alpha}\exp\left((1-\alpha)\int_{t_0}^{t}a(r)dr\right)+\\&(1-\alpha)\int_{t_0}^{t}b(s)\exp\left((1-\alpha)\int_{s}^{t}a(r)dr\right)ds \Bigg\}^{\frac{1}{1-\alpha}}
\end{align*}
\end{lemma}
 The proof of this lemma above can be found in p.360 of \cite{MPF}.

 Next we introduce a version of Gronwall's inequality with logarithmic nonlinearity.
 \begin{lemma}\label{6.2}
 	Let X, a, M, $c_1$, $c_2$ be nonnegative functions on $\mathbb{R}_+$, M be an increasing function and $M(0)>1$, and $c_1,\ c_2$ be locally integrable functions on $\mathbb{R}_+$. Assume that for any $t\geq 0$,\\
 	$$X(t)+a(t)\leq M(t)+\int_{0}^{t}c_1(s)X(s)ds+\int_{0}^{t}c_2(s)X(s)\log X(s)ds,$$
 	and the above integrals are finite. Then for any $t\geq 0$,
 	$$X(t)+a(t)\leq M(t)^{\exp(C_2(t))}\exp\left(\exp\left(C_2(t)\right)\int_{0}^{t}c_1(s)\exp(-C_2(s))ds\right)$$
 	where $C_2(t):=\int_{0}^{t}c_2(s)ds$.
  \end{lemma}
The proof can be seen in Lemma 7.2 of \cite{SZ}.
\vskip 0.4cm
\noindent {\bf Acknowledgement.} This work is partly supported by NSFC (No. 12131019, No. 11721101, No. 12001516)

\section*{Declarations}
The authors have no competing interests to declare that are relevant to the content of this article. Data sharing not applicable to this article as no datasets were generated or analysed during the current study.

\end{document}